\documentclass[12pt,a4paper]{article} 
\usepackage{amsmath,pstricks,latexsym,theorem,epsfig} 
\usepackage{amsfonts,amssymb,latexsym,epsfig,amsmath}
\usepackage{amsmath,latexsym,theorem} 
\newtheorem{Theorem}{Theorem}[section] 
\newtheorem{Definition}{Definition}[section] 
\newtheorem{Proposition}{Proposition}[section] 
\newtheorem{Lemma}{Lemma}[section] 
\newtheorem{Corollary}{Corollary}[section] 
 
\theorembodyfont{\rmfamily} 
\newtheorem{Remark}{Remark}[section]

\newcommand{\rec}[1]{{(\ref{#1})}} 
\def \N{I\!\!N} 

\def \R{I\!\!R} 
\def \C{\mathbb{C}}
\def\ds{\displaystyle}
 
\def \Z{Z\!\!\!Z} 
 
\def\11{1\!\!1}

 \def\res{\mathop{\hbox{\vrule height 7pt width .5pt 
depth 0pt\vrule height .5pt width 6pt depth 0pt}}\nolimits}

\def\rec#1{{(\ref{#1})}}

\newcommand{\ba}{\begin{array}} 
\newcommand{\ea}{\end{array}}

\begin{document} 
\title{\bf  3-Commutators Estimates and the Regularity of    $ {\mathbf {1/2}}$-Harmonic Maps into  Spheres} 
\author{ Francesca Da Lio\thanks{Department of Mathematics, ETH Z\"urich, R\"amistrasse 101, 8092 Z\"urich, Switzerland.} \thanks{Dipartimento di Matematica Pura ed Applicata, Universit\`a degli Studi di Padova. Via Trieste 63, 35121,Padova, Italy, e-mail: dalio@mah.unipd.it.} \and Tristan Riviere$^*$   }

\maketitle 
 
 \begin{abstract}
   We prove the regularity of weak $1/2-$harmonic maps  from the real line into a sphere. The key point in our result is first
   a formulation of the $1/2-$harmonic map equation in the form of  a non-local  linear Schr\"odinger type equation 
 with a  {\it 3-terms commutators } in the right-hand-side . We then establish a sharp estimate for these {\em 3-commutators}.
  \end{abstract}\par
  {\small {\bf Key words.} Harmonic maps, nonlinear elliptic PDE's, regularity of solutions, commutator estimates.}\par
 {\small { \bf  MSC 2000.}  58E20, 35J20, 35B65, 35J60, 35S99}
  \tableofcontents
\section{Introduction}
Since the early  50's the analysis of critical points to conformal invariant Lagrangians has raised a special interest, due to the important role they play in physics and geometry.\par
For a complete
overview on this topic we refer the reader to the introduction of \cite{Riv3}\,.
Here we recall  some classical  examples of  conformal invariant variational problems.  \par
The most elementary example of a $2$-dimensional conformal invariant Lagrangian is the Dirichlet Energy
\begin{equation}\label{lagr1}
E(u)=\int_{D}|\nabla u(x,y)|^2 dx dy\,,
\end{equation}
where $D\subseteq R^2$ is an open set and $u\colon D\to\R$, $\nabla u$ is the gradient of $u$\,.
We recall that a map $\phi\colon \C\to \C$ is conformal if it  satisfies
\begin{equation}\label{conf}
\left\{\begin{array}{l}
\displaystyle{|\frac{\partial\phi}{\partial x}|=|\frac{\partial\phi}{\partial y}|}
\\
\displaystyle{\langle \frac{\partial\phi}{\partial x},\frac{\partial\phi}{\partial y}\rangle=0}\\
{\rm det}\nabla\phi\ge 0 ~~{\mbox and} ~~\nabla\phi\ne 0
 \end{array} \right.
\end{equation}
Here   $\langle\cdot,\cdot\rangle$ denotes the
standard Euclidean inner product in $\R^n$\,.\par
For every $u\in W^{1,2}(D,\R)$ and every conformal map $\phi$, ${\rm deg}(\phi)=1,$ the following holds
$$
E(u)=E(u\circ\phi)=\int_{\phi^{-1}(D)}|(\nabla\circ\phi) u(x,y)|^2 dx dy\,.
$$
Critical points of this functional are the harmonic functions satisfying 
\begin{equation}\label{lapl}
\Delta u=0, ~~\mbox{in $ D$}\,.\end{equation}
We can extend $E$ to maps taking values in $\R^m$ as follows
\begin{equation}\label{lagr2}
E(u)=\int_{D}|\nabla u(x,y)|^2 dx dy=\int_{D}\sum_{i=1}^m|\nabla u_i(x,y)|^2 dx dy\,,
\end{equation}
where $u_i$ are the components of $u\,.$
The Lagrangian \rec{lagr2} is still conformally invariant and each component of its critical points satisfies the equation \rec{lapl}\,.\par
We can define the Lagrangian \rec{lagr2} also in the set of maps taking values in a 
compact submanifold ${\cal{N}}\subseteq\R^m$ without boundary. In this case   critical points $u\in W^{1,2}(D,{\cal{N}})$ of $E$   
  satisfy in a weak sense the equation
  $$
  -\Delta u\perp T_{u}{\cal{N}}\,,
  $$
where  $ T_\xi{\cal{N}}$ is the tangent plane a $\cal{N}$ at the point $\xi\in {\cal{N}}$, 
 or in a equivalent way 
  \begin{equation}\label{harm}
-\Delta u=A(u)(\nabla u,\nabla u):=A(u)(\partial_x u,\partial_x   u)+A(u)(\partial_y u,\partial_y   u),
\end{equation}
where $A(\xi)$ is the second fundamental form at the point $\xi\in\cal{N}$ (see for instance \cite{Hel}).
The equation \rec{harm} is called the {\em harmonic map equation} into ${\cal{N}}\,.$\par
In the case when $\cal{N}$ is an oriented hypersurface of $\R^m$ the harmonic map equation reads as
\begin{equation}\label{harmhyp}
-\Delta u=n\langle \nabla n,\nabla u\rangle\,,
\end{equation}
where $n$ denotes the composition of $u$ with the unit normal vector field $\nu$  to $\cal{N}\,.$\par
All the above examples belongs to the class of conformal invariant coercive Lagrangians whose 
corresponding Euler-lagrangian equation is of the form
\begin{equation}\label{EL}
-\Delta u=f(u,\nabla u)\,,
\end{equation}
where $f\colon\R^2\times(\R^m\otimes R^2)\to \R^m$ is a continuous function satisfying  for some positive constant 
$C$ 
$$
C^{-1}|p|^2\le f(\xi,p)|\le C|p|^2,~~~\forall~\xi,p\,.$$
One of the main issues related to equation \rec{EL} is the regularity of solutions $u\in W^{1,2}(D,{\cal{N}})$. We observe that equation \rec{EL} is critical in dimension 
$n=2$ for the $W^{1,2}$-norm. Indeed if we plug in the nonlinearity $f(u,\nabla u)$ the information that $u\in W^{1,2}(D,{\cal{N}})$, we get that
$\Delta u\in L^1(D)$ and thus $\nabla u\in L^{2,\infty}_{loc}(D)$ the weak $L^2$ space (see \cite{Ste}), which has the same homogeneity of $L^2$\,. Hence we are back in some sense to the initial situation. This shows that the equation is critical.\par
In general $W^{1,2}$ solutions to equations \rec{EL} are not smooth in dimension greater that $2$ (see
counter-example in \cite{Riv}).  We refer again the reader to \cite{Gia} for a more complete presentation of the
results concerning the regularity and compactness results for  equations \rec{EL}.\par
Here we are going to recall the approach introduced by F. H\'elein \cite{Hel} to prove the regularity of harmonic maps from a domain $D$ of $\R^2$ into the unit sphere $S^{m-1}$ of $\R^m$\,. In this case the 
Euler-Lagrange equation is
\begin{equation}
\label{harmsp}
-\Delta u=u|\nabla u|^2\,.
\end{equation}
It was observed by Shatah \cite{Sha} that  $u\in W^{1,2}(D,S^{m-1})$ is a solution of
\rec{harmsp} if and only if the following conservation law holds
\begin{equation}\label{cl}
{\rm div}(u_i\nabla u_j-u_j\nabla u_i)=0,
\end{equation}
for all $i,j\in\{1,\ldots,m\}\,.$\par
Using \rec{cl} and the fact that $|u|\equiv 1$ $\Longrightarrow$ $\sum_{j=1}^m u_j\nabla u_j=0$, H\'elein wrote the equation \rec{harmsp} in the form
\begin{equation}\label{harmsp2}
-\Delta u=\nabla^\perp B\cdot \nabla u,
\end{equation}
where $\nabla^\perp B=(\nabla^\perp B_{ij})$  with $\nabla^\perp B_{ij}=u_i\nabla u_j-u_j\nabla u_i\,,$ 
(for every
vector field  $v\colon \R^2\to \R^n$, $\nabla^\perp v$ denotes the $\pi/2$ rotation of the gradient $\nabla v$, namely  $\nabla^\perp v=(-\partial_y v,\partial_x v))\,.$\par
The r.h.s of \rec{harmsp2} can be written as a sum of jacobians:
$$
\nabla^\perp B_{ij}\nabla u_j=\partial_xu_j\partial_y B_{ij}-\partial_yu_j\partial_x B_{ij}\,.
$$
This particular structure permits to apply to the equation \rec{harmsp} the following result
\begin{Theorem}
\label{th-I.1}
\cite{W}
Let $D$ be a smooth  bounded domain  of ${\R}^2$. Let $a$ and $b$ be two measurable functions in $D$ whose gradients are in $L^2(D)$. Then there exists a unique solution $\varphi\in W^{1,2}(D)$ to
\begin{equation}\label{jac}
\left\{\begin{array}{ll}
-\Delta\varphi =\displaystyle{\frac{\partial a}{\partial x} \frac{\partial b}{\partial y} -\frac{\partial a}{\partial y} \frac{\partial b}{\partial x}},&~\mbox{in $D$}\\
\varphi=0 &~\mbox{on $\partial D\,.$}
\end{array}\right.
\end{equation}
Moreover there exists a constant $C>0$ independent of $a$ and $b$ such that
$$
||\varphi||_\infty+||\nabla\varphi||_{L^{2}}\le C||\nabla a||_{L^2}||\nabla b||_{L^2}\,.
$$
In particular $\varphi$ is a continuous in $D\,.$
\end{Theorem}
Theorem~\ref{th-I.1} applied to equation (\ref{harmsp2}) leads, modulo some standard localization argument in elliptic PDE, to an estimate of the form
\begin{equation}
\label{zz1}
\|\nabla u\|_{L^2(B_r(x_0))}\le C\ \|\nabla B\|_{L^2(B_r(x_0))}\ \|\nabla u\|_{L^2(B_r(x_0))}+Cr\ \|\nabla u\|_{L^2(\partial B_r(x_0))}
\end{equation}
for every $x_0\in D$ and $r>0$ such that $B_r(x_0)\subset D$. Assuming we are considering radii $r<r_0$ such that $ \max_{x_0\in D}C\ \|\nabla B\|_{L^2(B_{r}(x_0))}<1/2$, then
(\ref{zz1}) implies a Morrey estimate of the form
\begin{equation}
\label{zz2}
\sup_{x_0, r>0}\ r^{-\beta}\ \int_{B_r(x_0)}|\nabla u|^2\ dx<+\infty
\end{equation}
for some $\beta>0$ which itself implies the H\"older continuity of $u$ by standard embedding result (see \cite{Gia}). Finally a bootstrap argument implies that $u$ is in fact $C^\infty$ - and even analytic -
(see \cite{HW} and \cite{Mo}).

\medskip

In the present work we are interested in $1$ dimensional quadratic Lagrangians which are invariant under the trace of conformal maps that keep invariant the half space $\R^2_+$: the M\"oebius group.
\par
A typical example is the following
Lagrangian that we will call {\it $L-$energy} - $L$ stands for ''Line'' -
\begin{equation}\label{lagr3}
L(u)=\int_{\R}|\Delta^{1/4} u(x)|^2 dx\,,
\end{equation}
where $u\colon\R\to{\cal{N}}$, ${\cal{N}}$ is  a smooth   $k$-dimensional submanifold of  $\R^m$   which is at least $C^2$, compact and without boundary.  
We observe that the $L(u)$ in \rec{lagr3}  coincides with the semi-norm $||u||^2_{\dot{H}^{1/2}(\R)}$ (for the definition of $||\cdot ||_{\dot{H}^{1/2}(\R)}$ we refer to Section \ref{defnot})\,. Moreover a more
tractable way to look at this norm is given by the following identity
 \begin{eqnarray*}
\int_{\R}|\Delta^{1/4} u(x)|^2 dx=\inf\left\{\int_{\R^2_+}|\nabla \tilde u|^2 dx:~\tilde u\in W^{1,2}(\R^2,\R^m),~~\mbox{trace $\tilde u=u$}\right\}\,.
\end{eqnarray*}

The Lagragian $L$ extends to map $u$ in the following function space
 $$
 \dot{H}^{1/2}(\R,{\cal{N}})=\{u\in \dot{H}^{1/2}(\R,\R^m):~~u(x)\in {\cal{N}}, {\rm a.e},\}\,.
 $$
The operator $\Delta^{1/4} $ on ${\R}$ is defined by means of the the Fourier tranform as follows  $$\widehat{\Delta^{1/4} u}=|\xi|^{1/2}\hat u\,,$$
 (given a function f, $\hat f$ denotes the Fourier transform of $f$).\par

Denote $\pi_{\cal{N}}$ the orthogonal projection onto ${\cal{N}}$ which happens to be a $C^l$ map in a sufficiently small
neighborhood of ${\cal{N}}$ if ${\cal{N}}$ is assumed to be $C^{l+1}$. We now introduce the notion of $1/2$-harmonic map into a manifold.

 \begin{Definition}
 \label{df1}
 A map $u\in \dot{H}^{1/2}(\R,{\cal{N}})$   is called a weak $1/2$-harmonic map into $\cal{N}$ if for any $\phi\in C^\infty_0({\R},{\R}^m)$ there holds
 \[
 \frac{d}{dt}L(\pi_{\cal{N}}(u+t\phi))=0\quad.
 \]
 \hfill $\Box$
 \end{Definition}
 In short we say that a weak $1/2-$harmonic map is a {\it critical point of $L$ in  $\dot{H}^{1/2}(\R,{\cal{N}})$ for perturbations in the target}.

\medskip

 $1/2-$harmonic maps into the circle $S^1$   might appear for instance in the asymptotic of equations in phase-field theory for fractional reaction-diffusion such as
 \[
\epsilon^2\,\Delta^{1/2}u+u (1-|u|^2)=0
 \]
 where $u$ is a complex valued ''wave function''.
 
\medskip

In this paper we consider the case ${\cal{N}}=S^{m-1}$. We first write the Euler-Lagrange equation associated to $L$ in $\dot{H}^{1/2}(\R,S^{m-1})$
in the following way

\begin{Proposition}
\label{pr-I.1}
A map $u$ in $\dot{H}^{1/2}(\R,S^{m-1})$ is a weak $1/2$-harmonic map if and only if it satisfies the following Euler-Lagrange equation
\begin{equation}\label{eq1intr}
\Delta^{1/4} ( u\wedge \Delta^{1/4}u)=T(u\wedge,u)\,,
\end{equation}
where, in general for an arbitrary integer $n$, for every $Q\in \dot{H}^{1/2}(\R^n, {\cal{M}}_{\ell\times m}(\R))$ $\ell\ge 0$\footnote{ ${\cal{M}}_{\ell\times m}(\R)$ denotes, as usual, the space of $\ell\times m$ real matrices.} and $u\in\dot{H}^{1/2}(\R^n,\R^m)$ , $T$ is the operator defined by
\begin{equation}\label{opT}
T(Q,u):=\Delta^{1/4}(Q\Delta^{1/4}u)-Q\Delta^{1/2} u+\Delta^{1/4} u\Delta^{1/4} Q\,.\end{equation}
\hfill $\Box$
\end{Proposition}

The Euler Lagrange equation (\ref{eq1intr}) will often be completed by the following ''structure equation'' which a consequence of the fact that $u\in{S^{m-1}}$ almost everywhere :

\begin{Proposition}
\label{pr-I.2}
All maps in $\dot{H}^{1/2}(\R,S^{m-1})$ satisfy the following identity
\begin{equation}\label{eq2intr}
\Delta^{1/4} (u\cdot  \Delta^{1/4}u)=S(u\cdot,u)-{\cal{R}}(\Delta^{1/4} u \cdot {\cal{R}}\Delta^{1/4} u)\,.
\end{equation}
where,  in general for an arbitrary integer $n$, for every  $Q\in \dot{H}^{1/2}(\R^n, {\cal{M}}_{\ell\times m}(\R))$, $\ell\ge 0$ and $u\in\dot{H}^{1/2}(\R^n,\R^m)$, $S$ is the operator given by
\begin{equation}
\label{opS}
S(Q,u):=\Delta^{1/4}[Q\Delta^{1/4} u]-{\cal{R}} (Q\nabla u)+{\cal{R}}(\Delta ^{1/4} Q{\cal{R}}\Delta ^{1/4} u)\Delta^{1/4} u \cdot {\cal{R}}\Delta^{1/4} u
\end{equation}
and   ${\cal{R}}$ is  the Fourier multiplier of symbol $m(\xi)=i\frac{\xi}{|{\xi}|}\,$.
\hfill $\Box$
\end{Proposition}

\par 
 
In the present work we will first show that $\dot{H}^{1/2}$ solutions to the $1/2$-harmonic map equation (\ref{eq1intr}) are H\"older continuous. This regularity result 
will be a direct consequence of the following Morrey type estimate that we will establish :

\begin{equation}\label{holdcond}
\sup_{x_0\in{\R}, r>0}r^{-\beta}\int_{B_r(x_0)} |\Delta^{1/4}u|^2\ dx<+\infty
\end{equation}

To this purpose, in the spirit of what we have just presented regarding H\'elein's proof of the regularity of harmonic maps from a 2-dimensional domain
into a round sphere, we will take advantage of a ''gain of regularity'' in the r.h.s of the equations (\ref{eq1intr}) and (\ref{eq2intr}) 
where the different terms  $T(u\wedge, u)$, $S(u\cdot,u)$ and ${\cal{R}}(\Delta^{1/4} u \cdot {\cal{R}}\Delta^{1/4} u)$ play more or less the role
which was played by $\nabla^\perp B\cdot \nabla u$ in (\ref{harmsp2}). Precisely we will establish the following estimates : for every $u\in \dot{H}^{1/2}({\R},{\R}^m)$ and 
$ Q\in H^{1/2}(\R,{\cal{M}}_{\ell\times m}(\R))$ we have
\begin{equation}
\label{zz3}
\|T(Q,u)\|_{H^{-1/2}}\le C\ \|Q\|_{\dot{H}^{1/2}(\R)}\ \|u\|_{\dot{H}^{1/2}(\R)}\,,
\end{equation}
\begin{equation}
\label{zz4}
\|S(Q,u)\|_{H^{-1/2}}\le C\ \|Q\|_{\dot{H}^{1/2}(\R)}\ \|u\|_{\dot{H}^{1/2}(\R)}\,,
\end{equation}
and
\begin{equation}
\label{zz5}
\|{\cal{R}}(\Delta ^{1/4} u\cdot {\cal{R}}\Delta ^{1/4} u))\|_{\dot{H}^{-1/2}}\le C\  \|u\|_{\dot{H}^{1/2}(\R)}^2\,.
\end{equation}
Our denomination ''gain of regularity'' has been chosen in order to illustrate that, under our assumptions $u\in \dot{H}^{1/2}({\R},{\R}^m)$ and 
$Q\in \dot{H}^{1/2}(\R,{\cal{M}}_{\ell\times m}(\R))$ each term individually in $T$ and $S$ - like for instance $\Delta^{1/4}(Q\Delta^{1/4}u)$ or $Q\Delta^{1/2} u$ ... - are not in $H^{-1/2}$
but the special linear combination of them constituting $T$ and $S$ are in $H^{-1/2}$.
In a similar way, \underbar{in dimension 2},  $J(a,b):={\frac{\partial a}{\partial x} \frac{\partial b}{\partial y} -\frac{\partial a}{\partial y} \frac{\partial b}{\partial x}}$ satisfies, as a direct consequence of Wente's theorem
above,
\begin{equation}
\label{zz6}
\|J(a,b)\|_{\dot{H}^{-1}}\le C\ \|a\|_{\dot{H}^1}\ \|b\|_{\dot{H}^1}
\end{equation}
whereas, individually, the terms $\frac{\partial a}{\partial x} \frac{\partial b}{\partial y}$ and $\frac{\partial a}{\partial y} \frac{\partial b}{\partial x}$ are not in $H^{-1}$.

\medskip

The estimates (\ref{zz3}) and (\ref{zz4}) are in fact consequences of the following {\it 3-terms commutator} or simply {\it  3-commutator estimates} which are valid in arbitrary dimension $n$
and which represent two of the main results of the present paper. We recall that $BMO$ denotes the space of {\it Bounded Mean Oscillations} functions of John and Nirenberg (see for instance
\cite{Gra2})
\[
\|u\|_{BMO({\R}^n)}=\sup_{\{x_0\in {\R}^n\ ;\ r>0\}}\frac{1}{|B_r(x_0)|}\int_{B_r(x_0)}\left|u(x)-\frac{1}{|B_r(x_0)|}\int u(y)\ dy\right|\ dx\quad.
\]

   \begin{Theorem}\label{commestintr1} Let $n\in {\N}^\ast$ and
 let $u\in BMO(\R^n)$, $Q\in \dot{H}^{1/2}({\R}^n,{\cal{M}}_{\ell\times m}(\R))$ . Denote
 \[
 T(Q,u):=\Delta^{1/4}(Q\Delta^{1/4}u)-Q\Delta^{1/2} u+\Delta^{1/4} u\Delta^{1/4} Q\quad ,
 \]
 then $T(Q,u)\in {H}^{-1/2}$ and there exists $C>0$, depending only on $n,$ such that
\begin{equation}
\label{zz7}
||T(Q,u)||_{{{H}}^{-1/2}}\le C\ ||Q||_{\dot{H}^{1/2}(\R)}||u||_{BMO}\,.
\end{equation}
\hfill$\Box$
\end{Theorem}
  \begin{Theorem}\label{commestintr2} Let $n\in{\N}^\ast$ and 
 let $u\in BMO(\R^n)$, $Q\in \dot{H}^{1/2}({\R}^n,{\cal{M}}_{\ell\times m}(\R))$ . Denote
 \[
 S(Q,u):=\Delta^{1/4}[Q\Delta^{1/4} u]-{\cal{R}}  (Q\nabla u)+{\cal{R}}(\Delta ^{1/4} Q{\cal{R}}\Delta ^{1/4} u)
 \]
 where ${\cal{R}}$ is  the Fourier multiplier of symbol $m(\xi)=i\frac{\xi}{|{\xi}|}\,$. Then $S(Q,u)\in{H}^{-1/2}$ and there exists $C$ depending only on $n$ such that
\begin{equation}
\label{zz8}
||S(Q,u)||_{{{H}}^{-1/2}}\le C\ ||Q||_{\dot{H}^{1/2}(\R)}||u||_{BMO}\,.
\end{equation}
\hfill$\Box$
\end{Theorem}
The fact that Theorem~\ref{commestintr1} and Theorem~\ref{commestintr2} imply estimates (\ref{zz3}) and (\ref{zz4}) comes from the embedding $\dot{H}^{1/2}({\R})\subset BMO({\R})$.

The parallel between the structures $T$ and $S$ for $H^{1/2}$ in one hand and the jacobian structure $J$ for $H^1$ in the other hand can be pushed further as follows.
As a consequence of a result of R. Coifman, P.L. Lions, Y. Meyer and S. Semmes \cite{CLMS}, Wente estimate (\ref{zz6}) can be deduced from a
more general one : We denote, for any $i,j\in\{1\cdots n\}$, and $a,b\in \dot{H}^1({\R}^n)$,
\[
J_{ij}(a,b):=\frac{\partial a}{\partial x_i}\frac{\partial b}{\partial x_j}-\frac{\partial a}{\partial x_j}\frac{\partial b}{\partial x_i}\,,
\]
and denote $J(a,b):=(J_{ij}(a,b))_{ij=1\cdots n}$. With this notation the main result in \cite{CLMS} implies
\begin{equation}
\label{zz6bis}
\|J(a,b)\|_{\dot{H}^{-1}({\R}^n)}\le C\ \|a\|_{\dot{H}^1({\R}^n)}\ \|b\|_{BMO({\R}^n)}
\end{equation}
which is reminiscent to (\ref{zz7}) and (\ref{zz8}). Recall also that (\ref{zz6bis}) is a consequence of a commutator estimate by  
R. Coifman, R Rochberg and G. Weiss \cite{CRW}.

\medskip

The two theorems~\ref{commestintr1} and~\ref{commestintr1} will be the consequence of the two following ones
which are their ''dual versions''. Recall first that ${\mathcal H}^1({\R}^n)$ denotes the Hardy space of $L^1$
functions $f$ on ${\R}^n$satisfying 
\[
\int_{{\R}^n}\sup_{t\in {\R}}|\phi_t\ast f|(x)\ dx<+\infty\quad ,
\]
where $\phi_t(x):=t^{-n}\ \phi(t^{-1}x)$ and where $\phi$ is some function in the Schwartz space ${\mathcal S}({\R}^n)$ satisfying $\int_{{\R}^n}\phi(x)\ dx=1$. Recall the famous result by Fefferman 
saying that the dual space to ${\mathcal H}^1$ is $BMO$. 

 In one hand theorem~\ref{commestintr1} is the consequence of the following result.
\begin{Theorem}\label{comm1}
Let $u,Q\in \dot H^{1/2}(\R^n)$, denote
$$
R(Q,u)=\Delta^{1/4}(Q\Delta^{1/4} u)-\Delta^{1/2}(Qu)+\Delta^{1/4}((\Delta^{1/4} Q) u)\,.$$
then $R(Q,u)\in {\cal{H}}^{1}(\R^n)$ and
\begin{equation}\label{commest2}
||R(Q,u)||_{{\cal{H}}^{1}}\le C||Q||_{\dot{H}^{1/2}(\R)}||u||_{\dot{H}^{1/2}(\R)}\,.\end{equation}
\end{Theorem}
In the other hand theorem~\ref{commestintr2} is the consequence this next result.

\begin{Theorem}\label{comm3}
Let $u,Q\in H^{1/2}$ and $u\in BMO$. $$
\tilde S(Q,u)=\Delta^{1/4}(Q\Delta^{1/4} u)-\nabla(Q {\cal{R}}u)+{\cal{R}}\Delta^{1/4}(\Delta^{1/4} Q {\cal{R}}u)\,.$$
where  ${\cal{R}}$ is  the Fourier multiplier of symbol $m(\xi)=i\frac{\xi}{|{\xi}|}\,.$Then $\tilde S(Q,u)\in {\cal{H}}^1$ and 
\begin{equation}\label{commest3}
||\tilde S(Q,u)||_{{\cal{H}}^{1}}\le C||Q||_{\dot{H}^{1/2}(\R)}||u||_{\dot{H}^{1/2}(\R)}\,.\end{equation}
\hfill $\Box$

\end{Theorem}

\medskip

We now say few words on the proof of estimates~\ref{commest2} and~\ref{commest3}.
The compensations of the 3 different terms in $R(Q,u)$ will be clear from the Littlewood-Paley decomposition of the different products that we present in section 3. Denoting as usual $\Pi_1(fg)$ the high-low contribution - respectively from $f$ and $g$ - denoting $\Pi_2(fg)$ the low-high contribution and $\Pi_3(fg)$ the high-high contribution we shall need the following groupings

\begin{itemize}
\item{i)} For $\Pi_1(R(Q,u))$ we proceed to the following decomposition
$$
\Pi_1(R(Q,u))=\underbrace{\Pi_1(\Delta^{1/4}(Q\Delta^{1/4} u))}+\underbrace{\Pi_1(-\Delta^{1/2}(Qu)+\Delta^{1/4}((\Delta^{1/4} Q) u))}\,.$$
\item{ii)} For $\Pi_2(R(Q,u))$ we decompose as follows
$$
\Pi_2(R(Q,u))=\underbrace{\Pi_2(\Delta^{1/4}(Q\Delta^{1/4} u)-\Delta^{1/2}(Qu))}+\underbrace{\Pi_2(\Delta^{1/4}((\Delta^{1/4} Q) u))}\,.$$
\item{ii)} Finally, for $\Pi_3(R(Q,u))$ we decompose as follows
$$
\Pi_3(R(Q,u))=\underbrace{\Pi_3(\Delta^{1/4}(Q\Delta^{1/4} u))}-\underbrace{\Pi_3(\Delta^{1/2}(Qu))}+\underbrace{\Pi_3(\Delta^{1/4}((\Delta^{1/4} Q) u))}\,.$$
\end{itemize}


Finally, injecting the Morrey estimate (\ref{holdcond}) in equations (\ref{eq1intr}) and (\ref{eq2intr}), a classical ''elliptic type'' bootstrap argument leads to the following result  
(see \cite{DR} for the details of this argument).
\begin{Theorem}
\label{th-I.3}
Let $u$ be a weak $1/2$-harmonic map in $\dot{H}^{1/2}({\R},S^{m-1})$. Then it belongs to $H^{s}_{loc}(\R,{S^{m-1}})$ for every $s\in{\R}$ and thus it is $C^\infty$ .
\hfill $\Box$
\end{Theorem}

\medskip

The paper is organized as follows.

\begin{itemize}

\item[-] In Section \ref{defnot} we give some preliminary definitions and notations.

\item[-] In Section \ref{commsection} we prove the {\it 3-commutator estimates} Theorems \ref{commestintr1} and \ref{commestintr2}. 
 
\item[-] In section 4 we study geometric localization properties of the $\dot{H}^{1/2}-$ norm on the real line
 for $\dot{H}^{1/2}-$functions in general
 
\item[-] In Section \ref{regulsection} we prove some $L-$energy decrease control on dyadic annuli for general solutions to 
 some linear non-local systems of equations that will include the systems (\ref{eq1intr}) and (\ref{eq2intr}).
 
 \item[-] in Section \ref{application} we derive the Euler-Lagrange equation (\ref{eq1intr}) associated to the Lagrangian 
\rec{lagr3} - proposition~\ref{pr-I.1}. We then prove proposition~\ref{pr-I.2}. We finally use the results of the previous section in order to deduce the Morrey type estimate (\ref{holdcond}) for $1/2-$harmonic maps into a sphere\,.


\end{itemize}


\section{Notations and Definitions}\label{defnot}

In this Section we introduce some  notations and definitions we are going to use in the sequel.
\par

For $n\ge 1$,  we denote respectively by ${\cal{S}}(\R^n)$ and  ${\cal{S}}^{\prime}(\R^n)$ the spaces of Schwartz functions and tempered distributions. Moreover given a function $v$ we will denote either by
  $\hat v$ or by  ${\mathcal{F}}[v]$ the Fourier Transform of  $v$ :
  $$\hat v(\xi)={\mathcal{F}}[v](\xi)=\int_{\R^n}v(x)e^{-i \langle \xi, x\rangle }\,dx\,.$$
  Throughout the paper we use the convention that $x,y$ denote variables in the space and
  $\xi,\zeta$ variables in the phase\,.
  
  We recall the definition of fractional Sobolev space (see for instance \cite{T3}).\par
  \begin{Definition}\label{fracsob} 
  For a real $s\ge 0$, 
  $$H^{s}(\R^n)=\{v\in L^2(\R^n):~~|\xi|^s{\cal{F}}[v]\in L^2(\R^n)\}.$$
  For a real $s<0$,
   $$H^{s}(\R^n)=\{v\in {\cal{S}}^{\prime} (\R^n):~~(1+|\xi|^2)^s{\cal{F}}[v]\in L^2(\R^n)\}.$$
   \end{Definition}
   It is known that $H^{-s}(\R^n)$ is the dual of $H^{s}(\R^n)\,.$\par
   Another characterization of $H^{s}(\R^n)$, with $0<s<1$, which does not use the Fourier transform is the following, (see for instance \cite{T3}).
   \begin{Lemma}
   For $0<s<1$, $u\in H^{s}(\R^n)$ is equivalent to $u\in L^2(\R^n)$ and
   $$
 \left(\int_{\R^n}\int_{\R^n} \left(\frac{(u(x)-u(y))^2}{|x-y|^{n+2s}}\right)dx dy\right)^{1/2}<+\infty\,.$$
 \end{Lemma}
  For $s>0$ we set
  $$||u||_{H^s(\R^n)}=||u||_{L^2(\R^n)}+|||\xi|^s{\cal{F}}[v]||_{L^2(\R^n)}\,,$$
  and
  $$||u||_{\dot{H}^s(\R^n)}=|||\xi|^s{\cal{F}}[v]||_{L^2(\R^n)}\,.$$
  
 For an open set $\Omega\subset\R^n$,   $H^s(\Omega)$ is the space of the restrictions of functions from $H^{s}(\R^n)$ and
 $$
 ||u||_{\dot H^{s}(\Omega)}=\inf\{||U||_{\dot H^{s}(\R^n)},~~U=u ~\mbox{on}~\Omega\}$$
 In the case of  $0<s<1$ then 
  $f\in H^{s}(\Omega)$ if and only if 
 $f\in L^2(\Omega)$ and  
 $$
 \left(\int_\Omega\int_\Omega \left(\frac{(u(x)-u(y))^2}{|x-y|^{n+2s}}\right)dx dy\right)^{1/2}<+\infty\,.$$
 Moreover
 $$
  ||u||_{\dot H^{s}(\Omega)}\simeq \left(\int_\Omega\int_\Omega \left(\frac{(u(x)-u(y))^2}{|x-y|^{n+2s}}\right)dx dy\right)^{1/2}<+\infty\,,$$
  see for instance \cite{T3}\,.
  
  Finally for a submanifold  ${\cal{N}}$ of $\R^m$  we can define   
 $$H^{s}(\R,{\cal{N}})=\{u\in H^{s}(\R,\R^m):~~u(x)\in {\cal{N}}, {\rm a.e.}\}\,.
 $$
  We introduce  the so-called Littlewood-Paley or dyadic  decomposition of unity. Such a decomposition can be obtained as follows \,.
 Let $\phi(\xi)$ be a radial Schwartz function supported on $\{\xi:~|\xi|\le 2\}$, which is
 equal to $1$ on $\{\xi: ~|\xi|\le 1\}$\,.
 Let $\psi(\xi)$ be the function
 $
 \psi(\xi):=\phi(\xi)-\phi(2\xi)\,.$
 $\psi$ is a bump function supported on the annulus $\{\xi:~1/2\le |\xi|\le 2\}\,.$\par
 We put $\psi_0=\phi$, $\psi_j(\xi)=\psi(2^{-j}\xi)$ for $j\ne 0 \,.$ The functions $\psi_j$, for $j\in\Z$, are supported on  $\{\xi:~2^{j-1}\le |\xi|\le 2^{j+1}\}\,.$ Moreover $\sum_{j\in\Z}\psi_j(x)=1\,.$\par 
 We then set $\phi_j(\xi):=\sum_{k=-\infty}^j\psi_k(\xi)\,.$ 
 The function $\phi_j$ is supported on  $\{\xi, ~|\xi|\le 2^{j+1}\}$.\par

  We recall the definition of the homogeneous Besov spaces $\dot{B}_{p,q}^s(\R^n)$  and  homogeneous Triebel-Lizorkin spaces $\dot{F}_{pq}^s(\R^n)$ in terms of the above dyadic decomposition.
  \begin{Definition}
  Let $s\in\R$,  $0< p,q\le\infty\,.$ For $f\in{\cal{S}}^\prime (\R^n)$ we set 
   \begin{equation}\left.\begin{array}{ll}
  ||f||_{\dot{B}_{p,q}^s(\R^n)}=\left(\sum_{j=-\infty}^{\infty}2^{jsq}||{\cal{F}}^{-1}[\psi_j {\cal{F}}[f]]||_{L^{p}(\R^n)}^q\right)^{1/q}&~~\mbox{if $q<\infty$}\\
    ||f||_{\dot{B}_{p,q}^s(\R^n)}=\sup_{j\in \N} 2^{js}||{\cal{F}}^{-1}[\psi_j {\cal{F}}[f]]||_{L^{p}(\R^n)}&~~\mbox{if $q=\infty$}\end{array}\right.
    \end{equation}
    When $p,q<\infty$ we also set
   $$ ||f||_{\dot{F}_{p,q}^s(\R^n)}=||\left(\sum_{j=-\infty}^{\infty}2^{jsq}|{\cal{F}}^{-1}[\psi_j {\cal{F}}[f]]|^q\right)^{1/q}||_{L^p}\,.$$
    \end{Definition}
    The space of all tempered distributions $f$ for which the quantity $||f||_{\dot{B}_{p,q}^s(\R^n)}$ is finite is called the homogeneous Besov space with indices 
    $s,p,q$ and it is denoted by $\dot{B}_{p,q}^s(\R^n)$. The space of all tempered distributions $f$ for which the quantity $||f||_{\dot F_{p,q}^s(\R^n)}$ is finite is called the homogeneous
    Triebel-Lizorkin space with indices 
    $s,p,q$ and it is denote by $\dot F_{p,q}^s(\R^n)\,.$ 
    It is known that $\dot{H}^s(\R^n)=\dot{B}^s_{2,2}(\R^n)=\dot{F}_{2,2}^s(\R^n)$\,.  
    
  Finally we denote    ${\cal{H}}^1(\R^n)$  the homogeneous Hardy Space in $\R^n$. It is known that ${\cal{H}}^1(\R^n)\simeq F^0_{2,1}$ thus we have  
     $$||f||_{{\cal{H}}^{1}(\R^n)}\simeq \int_{\R}(\sum_j |{\cal{F}}^{-1}[\psi_j {\cal{F}}[f]]|^2)^{1/2}dx\,.$$
  \par
  
  We recall that in dimension $n=1$,  the space  $\dot{H}^{1/2}(\R)$ is continuously embedded in the Besov space $ \dot{B}^0_{\infty,\infty}(\R)$.
 More precisely we have 
 \begin{equation} 
 \dot{H}^{1/2}(\R)\hookrightarrow BMO(\R)\hookrightarrow \dot{B}^0_{\infty,\infty}(\R)\,,\end{equation}
 where  (see for instance page 31 in  \cite{RS}, page 129 in \cite{Tr1}).
 \par 
 The $s$-fractional  Laplacian of a function  $u\colon\R^n\to\R$ is defined as a pseudo differential operator of symbol $|\xi|^{2s}$ :
 \begin{equation}\label{fract}
 \widehat {\Delta^{s}u}(\xi)=|\xi|^{2s} \hat u(\xi)\,.
 \end{equation}
 In the case where $s=1/2 $, we can write ${\Delta}^{1/2}u= -{\cal{R}}(\nabla u)$ where ${\cal{R}} $ is  Fourier multiplier of symbol
 $\displaystyle{\frac{i}{|\xi|}\sum_{k=1}^n\xi_k}:$
 $$
 \widehat {{\cal{R}}X}(\xi)=\frac{1}{|\xi|}\sum_{k=1}^n i\xi_k\hat {X_k}(\xi)$$
 for every $X\colon \R^n\to\R^n$\,, namely $ {\cal{R}}={\Delta}^{-1/2}{\rm div}\,.$\par
   
  We denote by $B_r(\bar x)$ the ball of radius $r$ and centered at $\bar x$. If $\bar x=0$ we simply write
  $B_r$\,. If $x,y\in\R^n,$ $x\cdot y$ denote the scalar product between $x,y$\,.
   
  For every function $f\colon\R^n\to\R$ we denote by $M(f)$ the maximal function of $f$, namely
 \begin{equation}\label{maxf}
 M(f)=\sup_{r>0,x\in\R^n}|B(x,r)|^{-1}\int_{B(x,r)}|f(y)|dy\,.\end{equation}


\section{3-Commutator Estimates : Proof of theorem~\ref{commestintr1} and theorem~\ref{commestintr2}.}\label{commsection}

In this Section we prove Theorems \ref{commestintr1} and \ref{commestintr2}\,.\par
We consider the dyadic decomposition introduced in  Section \ref{defnot}\,.  
   For every $j\in \N$ and $f\in{\cal{S}}^{\prime}(\R^n)$  we set  
   \begin{eqnarray*}
  f^j={\cal{F}}^{-1}[\phi_j {\cal F}[f]],~~~
 f_j={\cal{F}}^{-1}[\psi_j {\cal F}[f]]\,.
 \end{eqnarray*}
 We have $f^j=\sum_{k=0}^{j} f_k$  and $f=\sum_{k=0}^{+\infty}f_k$ (where the convergence is in ${\cal{S}}^\prime(\R^n)$)\,.\par
 Let $f,g\in {\cal{S}}^\prime(\R)$. Suppose that $fg$ exists in ${\cal{S}}^\prime(\R^n)$. Then we split  the product in the following way
 $$
 fg=\Pi_1(fg)+\Pi_2(fg)+\Pi_3 (fg),$$
 where
 \begin{eqnarray*}
 \Pi_1(fg)&=& \sum_{-\infty}^{+\infty} f_j\sum_{k\le j-4} g_k= \sum_{-\infty}^{+\infty} f_j g^{j-4}\,;\\
 \Pi_2(fg)&=& \sum_{-\infty}^{+\infty} f_j\sum_{k\ge j+4} g_k \sum_{-\infty}^{+\infty} g_j f^{j-4}\,;\\
 \Pi_3(fg)&=& \sum_{-\infty}^{+\infty}  f_j\sum_{|k-j|< 4} g_k\,.\
 \end{eqnarray*}
 We observe that   for every $j$ we have 
 $$\mbox{supp${\cal{F}}[f^{j-4}g_j]\subset \{2^{j-2}\le |\xi|\le 2^{j+2}\}$};$$
 $$\mbox{supp${\cal{F}}[\sum_{k=j-3}^{j+3}f_jg_k]\subset \{|\xi|\le 2^{j+5}\}$}\,.$$
 The following Lemma will be often  used   in the sequel.
 \begin{Lemma}
 For every $f\in {\cal{S}}^{\prime}$ we have
 $$ \sup_{j\in Z}|f^j|\le M(f)\,.$$
 \end{Lemma}
 {\bf Proof.}
 We have
 \begin{eqnarray*}
 && 
 f^j={\cal{F}}^{-1}[\phi_j]\star f=2^j\int_R {\cal{F}}^{-1}[\phi](2^j(x-y))f(y)dy\\
 && =\int_{\R} {\cal{F}}^{-1}[\phi](z)f(x-2^{-j}z)dz\\
 &&=\sum_{k=-\infty}^{+\infty}\int_{B_{2^k}\setminus B_{2^{k-1}}}{\cal{F}}^{-1}[\phi](z)f(x-2^{-j}z)dz\\
 &&\le \sum_{k=-\infty}^{+\infty}\max_{B_{2^k}\setminus B_{2^{k-1}}}|{\cal{F}}^{-1}[\phi](z)|
 \int_{B_{2^k}\setminus B_{2^{k-1}}}|f(x-2^{-j}z)|dz\\
 &&\le  \sum_{k=-\infty}^{+\infty}\max_{B_{2^k}\setminus B_{2^{k-1}}} 2^{k}
 |{\cal{F}}^{-1}[\phi](z)|
2^{j-k} \int_{B(x,2^{k-j})\setminus B(x,2^{k-1-j})}|f(z)|dz\\
&&\le M(f)\sum_{k=-\infty}^{+\infty}\max_{B_{2^k}\setminus B_{2^{k-1}}}2^{k}|{\cal{F}}^{-1}[\phi](z)|\le C M(f)\,.
\end{eqnarray*}
In the last inequality we use the fact ${\cal{F}}^{-1}[\phi]$ is in ${\cal{S}}(\R^n)$ and thus
$$\sum_{k=-\infty}^{+\infty}\max_{B_{2^k}\setminus B_{2^{k-1}}}2^{k}|{\cal{F}}^{-1}[\phi](z)|\le 2\int_{\R}
|{\cal{F}}^{-1}[\phi](z)|d\xi,.$$

We can now start the proof of one of the main result in the paper.\par\medskip
{\bf Proof of theorem~\ref{comm1}.}\par
We are going to estimate $\Pi_1(R(Q,u)), \Pi_2(R(Q,u)) $ and $ \Pi_3(R(Q,u))\,.$\par
{$\bullet$}  Estimate of $||\Pi_1(\Delta^{1/4}(Q\Delta^{1/4} u)||_{{\cal{H}}^1}$\,.\par
\begin{eqnarray}\label{commestpi1}
&&||\Pi_1(\Delta^{1/4}(Q\Delta^{1/4} u)||_{{\cal{H}}^1}=\int_{\R^n}\left(\sum_{j=-\infty}^{\infty} 2^j Q^2_j(\Delta^{1/4} u^{j-4}\right)^2)^{1/2} dx\\
&&\le \int_{\R^n} \sup_{j}|\Delta^{1/4} u^{j-4}|\left(\sum_j 2^j Q_j^2\right)^{1/2} dx\nonumber\\
&&\le \left(\int_{\R^n} (M(\Delta^{1/4} u))^2 dx\right)^{1/2}\left(\int_R\sum_j 2^j Q_j^2 dx\right)^{1/2}\nonumber\\
&&\le C||u||_{\dot{H}^{1/2}}||Q||_{\dot{H}^{1/2}}\,.\nonumber
\end{eqnarray}

$\bullet$ Estimate of $\Pi_1(\Delta^{1/4}(\Delta^{1/4} Q u)-\Delta^{1/2}(Qu))$. We show that it is in
$\dot{B}^0_{1,1}$ (${\cal{H}}^1\hookrightarrow \dot{B}^0_{1,1}$)\,.
To this purpose
we use the ``commutator structure of the above term"\,.

 \begin{eqnarray}\label{p2com2}
&&||\Pi_1(\Delta^{1/4}(\Delta^{1/4} Q) u-\Delta^{1/2}(Qu))||_{\dot{B}^0_{1,1}}\\
&&=\sup _{||h||_{\dot{B}^0_{\infty,\infty}}\le 1} \int_{\R^n}\sum_j\sum _{|t-j|\le 3}[\Delta^{1/4}(u^{j-4}\Delta^{1/4} Q_j)-\Delta^{1/2}(u^{j-4}Q_j)]h_tdx\nonumber\\
&& =
\sup _{||h||_{\dot{B}^0_{\infty,\infty}}\le 1} \int_{\R^n}\sum_j\sum _{|t-j|\le 3} {\cal{F}}[u^{j-4}]
{\cal{F}}[\Delta ^{1/4}Q_j \Delta^{1/4} h_t-Q_j\Delta^{1/2}h_t]d\xi\nonumber\\ 
 && =\sup _{||h||_{\dot{B}^0_{\infty,\infty}}\le 1}\int_{\R^n}\sum_j\sum _{|t-j|\le 3} {\cal{F}}[u^{j-4}](\xi)\nonumber\\
&&~~~ \left(\int_{\R^n} {\cal{F}}[Q_j] (\zeta){\cal{F}}[\Delta^{1/4} h_t](\xi-\zeta)( |\zeta|^{1/2}-|\xi-\zeta|^{1/2} ) d\zeta\right) d\xi\,.\nonumber
\end{eqnarray}

Now we observe that in \rec{p2com2} we have $|\xi|\le 2^{j-3}$ and $2^{j-2}\le |\zeta|\le 2^{j+2}$.
Thus $|\displaystyle\frac{\xi}{\zeta}|\le \frac{1}{2}\,.$ 
\par
Hence
\begin{eqnarray}\label{estkern}
| |\zeta|^{1/2}-|\xi-\zeta|^{1/2} &=&|\zeta|^{1/2}[1-|1-\frac{\xi}{\zeta}|^{1/2}]\\
&=&|\zeta|^{1/2}\frac{\xi}{\zeta}[1+|1-\frac{\xi}{\zeta}|^{1/2}]^{-1}\nonumber\\
&=&|\zeta|^{1/2}\sum_{k=-\infty}^\infty\frac{c_k}{k!}(\frac{ \xi}{\zeta})^{k+1}\,.\nonumber
\end{eqnarray}
   We    introduce the following notation: for every $k\ge 0$  and $g\in{\cal{S}}^{\prime}$ we set
$$S_k g={\cal{F}}^{-1}[\xi^{-(k+1)}|\xi|^{1/2} {\cal{F}} g].$$
We note that  if $h\in B^s_{\infty,\infty}$ then
$S_k h\in \dot B^{s+1/2+k}_{\infty,\infty}$ and if $h\in \dot H^{s}$ then $S_k h\in \dot H^{s+1/2+k}\,.$\par
Finally  if $Q\in \dot H^{1/2}$ then $\nabla^{k+1}(Q)\in \dot H^{-k-1/2}\,.$\par
 We continue the estimate \rec{p2com2}\,.
 \begin{eqnarray*}
&&\sup _{||h||_{\dot{B}^0_{\infty,\infty}}\le 1}\int_{\R^n}\sum_j\sum _{|t-j|\le 3} {\cal{F}}[u^{j-4}](\xi)\\
&&~~~( \int_{\R^n}  {\cal{F}}[  Q_j] (\zeta){\cal{F}}[\Delta^{1/4} h_t](\xi-\zeta) (|\xi-\zeta|^{1/2}-(|\zeta|^{1/2}) d\zeta)  d\xi \\
&& =\sup _{||h||_{\dot{B}^0_{\infty,\infty}}\le 1}\int_{\R^n}\sum_j\sum _{|t-j|\le 3} {\cal{F}}[u^{j-4}]( \xi)\\
&&~~~\left[\int_{\R^n}|\zeta|^{1/2} {\cal{F}}[ Q_j] (\zeta){\cal{F}}[\Delta^{1/4}h_t](\xi-\zeta)\sum_{\ell=0}^\infty\frac{c_{\ell}}{\ell!}(\frac{\xi}{\eta})^{\ell+1} d\zeta\right]d\xi\\&&
\le C\sup _{||h||_{\dot{B}^0_{\infty,\infty}}\le 1} 
\sum_{\ell=0}^\infty\frac{c_{\ell}}{\ell!}\int_{\R^n}
\sum_j\sum _{|t-j|\le 3}
(i)^{-(\ell+1)}{\cal{F}}[\nabla^{\ell+1}u^{j-4}]{\cal{F}}[S_{\ell} Q_j\Delta^{1/4}h_t)]( \xi) d \xi\nonumber\\
&&
\le C
\sup _{||h||_{\dot{B}^0_{\infty,\infty}}\le 1} ||h||_{\dot{B}^0_{\infty,\infty}}
\sum_{\ell=0}^\infty\frac{c_{\ell}}{\ell!}\int_{\R^n}
\sum_j 2^{j/2}
|\nabla^{\ell+1}u^{j-4}||S_{\ell} Q_j|d x\nonumber\\
&& 
\le
\sum_{\ell=0}^\infty\frac{c_{\ell}}{\ell!}\int_{\R^n}
\sum_j |2^{-(k+1/2)j}\nabla^{\ell+1}u^{j-4}||2^{(k+1)j}S_{\ell} Q_j|dx\nonumber\\
&& 
\le  C \sum_{\ell=0}^\infty\frac{c_{\ell}}{\ell!}(\int_{\R^n}\sum_j 2^{-2(\ell+1/2)j}|\nabla^{\ell+1} u^{j-4}|^2 dx)^{1/2}(\int_{\R^n}\sum_j 2^{2(\ell+1)j}|S_{\ell} Q_j|^2 dx)^{1/2}\nonumber\\\
&&\mbox{by Plancherel Theorem}\\&&
=
C  \sum_{\ell=0}^\infty\frac{c_{\ell}}{\ell!}(\int_{\R^n}\sum_j 2^{-2(\ell+1/2)j}|\xi|^{2\ell}|{\cal{F}}[\nabla u^{j-4}]|^2 d\xi)^{1/2}(\int_{\R^n}\sum_j 2^{2(\ell+1)j}|\xi|^{-2(\ell+1/2)}|{\cal{F}}[Q_j]|^2 d\xi)^{1/2}
\nonumber\\\
&&
\le  
C  \sum_{\ell=0}^\infty\frac{c_{\ell}}{\ell!}2^{-3\ell}(\int_{\R^n}\sum_j 2^{-j}|{\cal{F}}[\nabla u^{j-4}]|^2 d\xi)^{1/2}(\int_{\R^n}\sum_j 2^{j}|{\cal{F}}[Q_j]|^2 d\xi)^{1/2}
\nonumber\\
&&\le C \sum_{\ell=0}^\infty\frac{c_{\ell}}{\ell!}2^{-3\ell} ||Q||_{\dot{H}^{1/2}}||u||_{\dot{H}^{1/2}}\,.\nonumber
\end{eqnarray*}

Above we also use the fact that for every vector field $X$ we have
\begin{eqnarray}\label{equiv}
&&\int_{\R^n}\sum_{j=-\infty}^{+\infty} 2^{-j}(X^{j-4})^2 dx
=\int_{\R^n} \sum_{k,\ell}X_kX_{\ell}\sum_{j-4\ge k,j-4\ge \ell}2^{-j}dx\nonumber\\ &&
\lesssim\int_{\R^n}\sum_{j=-\infty}^{+\infty} 2^{-j}(X_j)^2dx\,.\end{eqnarray}
 
 The estimate of $||\Pi_2(\Delta^{1/4}(Q\Delta^{1/4} u)-\Delta^{1/2}(Qu))||_{\dot{B}^0_{1,1}}$ is analogous to \rec{p2com2}\,.\par
{$\bullet$}  Estimate of $||\Pi_2(\Delta^{1/4}(\Delta^{1/4} Q u )||_{{\cal{H}}^1}$. It is as in \rec{commestpi1}\,.\par
{$\bullet$}  Estimate of $||\Pi_3(\Delta^{1/2}(Qu)||_{{\cal{H}}^1}$.\par We show that it is indeed in the smaller space  $\dot{B}^0_{1,1}$ ( we have $\dot{B}^0_{1,1}\hookrightarrow {{\cal{H}}^1}$)\,.
 We first observe that if $h\in \dot{B}^0_{\infty,\infty}$ then $\Delta^{1/2} h\in \dot{B}^{-1}_{\infty,\infty}$ and 
\begin{eqnarray*}
&&\Delta^{1/2} h^{j-6}=\sum_{k=-\infty}^{j}\Delta^{1/2} h_{k}\le \sup_{k\in\N} |2^{-k}\Delta^{1/2} h_k| \sum_{k=-\infty}^{j} 2^{k}\le C2^{j}||h||_{\dot{B}^{0}_{\infty,\infty}}\,.\\
 \end{eqnarray*}
\begin{eqnarray}\label{commestpi3}
&&
||\Pi_3(\Delta^{1/2}(Qu)||_{\dot{B}^0_{1,1}}=\sup _{||h||_{\dot{B}^0_{\infty,\infty}}\le 1}\int_{\R^n}\sum_j\sum_{|k-j|\le 3}\Delta^{1/2}(Q_j u_k)h\\
&&=\sup _{||h||_{\dot{B}^0_{\infty,\infty}}\le 1}\int_{\R^n}\sum_j\sum_{|k-j|\le 3}\Delta^{1/2}(Q_ju_k)
\left[ h^{j-6}+\sum_{t= j-5}^{j+6} h_t\right]dx\nonumber\\
 && =\sup _{||h||_{\dot{B}^0_{\infty,\infty}}\le 1}\int_{\R^n}\sum_j\sum_{|k-j|\le 3} (Q_ju_k)
\left[\Delta^{1/2} h^{j-6}+\sum_{t= j-5}^{j+6}\Delta^{1/2} h_t\right]dx\nonumber\\
&&\le C \sup _{||h||_{\dot{B}^0_{\infty,\infty}}\le 1}||h||_{\dot{B}^0_{\infty,\infty}}\int_{\R^n}\sum_j\sum_{|k-j|\le 3} 2^{j} |Q_j u_k| dx\nonumber\\
&&\le C\left(\int_{\R^n}\sum_j  2^{j} Q^2_j dx\right)^{1/2}\left(\int_{\R^n}\sum_j 2^{j} u^2_j dx\right)^{1/2}\nonumber\\
&&
\le C||Q||_{\dot{H}^{1/2}}||u||_{\dot{H}^{1/2}}\,.\nonumber
\end{eqnarray}
 
$\bullet$ Estimate of $\Pi_3(\Delta^{1/4}(Q\Delta^{1/4} u))$\,.\par
We show that it is in $\dot{B}^0_{1,1}$. \par
We observe that if $h\in \dot{B}^0_{\infty,\infty}$ then $\Delta^{1/4} h\in  B^{-1/2}_{\infty,\infty}$ and thus

\begin{eqnarray}\label{p3comm2}
&&||\Pi_3(\Delta^{1/4}(Q,\Delta^{1/4} u))||_{\dot{B}^0_{1,1}}=\sup _{||h||_{\dot{B}^0_{\infty,\infty}}\le 1}\int_{\R^n}\sum_j\sum_{|k-j|\le 3}\Delta^{1/4}(Q_j\Delta^{1/4} u_k)h\nonumber\\
&&=\sup _{||h||_{\dot{B}^0_{\infty,\infty}}\le 1}\int_{\R^n}\sum_j\sum_{|k-j|\le 3}\Delta^{1/4}(Q_j\Delta^{1/4} u_k)
\left[\Delta^{1/4} h^{j-6}+\sum_{t= j-5}^{j+6}\Delta^{1/4}h_t\right]dx\nonumber
 \\
&&\le C \sup _{||h||_{\dot{B}^0_{\infty,\infty}}\le 1}||h||_{\dot{B}^0_{\infty,\infty}}\int_{\R^n}\sum_j\sum_{|k-j|\le 3} 2^{j/2} |Q_j\Delta^{1/4} u_k| dx\\
&&\le C(\int_{\R^n}\sum_j  2^{j} Q^2_j dx)^{1/2}(\int_{\R^n}\sum_j (\Delta^{1/4} u_j)^2dx)^{1/2}\nonumber\\
&&
\le C||Q||_{\dot{H}^{1/2}}||u||_{\dot{H}^{1/2}}\,.\nonumber
\end{eqnarray}
 
The estimate of    $\Pi_3(\Delta^{1/4}(\Delta^{1/4} Q u))$ is analogous to \rec{p3comm2}\,.
~\hfill$\Box$\par\medskip
 From Theorem \ref{comm1} and the duality between $BMO$ and ${\cal{H}}^1$ we get Theorem \ref{commestintr1}.
 \par
 {\bf Proof of Theorem \ref{commestintr1}\,.}  
  
For all $h,Q\in \dot H^{1/2}$ and $u\in BMO$ we have
\begin{eqnarray*}
&& \int_{\R^n}[(\Delta^{1/4}(Q\Delta^{1/4}u)-Q\Delta^{1/2}u+\Delta^{1/4}Q\Delta^{1/4}u]h dx\\
&&  =\int_{\R^n} [(\Delta^{1/4}(Q\Delta^{1/4}h)-\Delta^{1/2}(Q h)+\Delta^{1/4}(h\Delta^{1/4}Q]u dx\\
&&~~~\mbox{by Theorem \rec{comm1}}\\
&&\le
C||u||_{BMO}||R(Q,h)||_{{\cal{H}}^1}\le C||u||_{BMO}||Q||_{\dot{H}^{1/2}}||h||_{\dot{H}^{1/2}}\,.
\end{eqnarray*}
Hence
\begin{eqnarray*}
&&
||T(Q,u)||_{{{\dot H}}^{-1/2}}=\sup_{||h||_{\dot{H}^{1/2}}\le 1}\int_{\R^n}T(Q,u) h dx\le C||u||_{BMO}||Q||_{\dot{H}^{1/2}}\,. 
\end{eqnarray*}
~~~\hfill$\Box$\par

{\bf Proof of theorem~\ref{comm3}.}
We observe that ${\cal{R}}$ is a Fourier multiplier of order zero thus ${\cal{R}}\colon{H^{-1/2}}\to {H^{-1/2}}$, ${\cal{R}}\colon{\cal{H}}^1\to {\cal{H}}^1$, and
 ${\cal{R}}\colon \dot{B}^0_{1,1}\to \dot{B}^0_{1,1}$
 (see \cite{Tay} and \cite{SiTri})\,.
  
The estimates are very similar to ones in Theorem \ref{comm1}, thus we will make   only the following one.\par

 $\bullet$ Estimate of $\Pi_1({\cal{R}}\Delta^{1/4}(\Delta^{1/4}Q {\cal{R}}u)-\nabla(Q{\cal{R}} u))$.\par
We observe that $\nabla u=\Delta^{1/4} {\cal{R}}\Delta^{1/4} u$
 \begin{eqnarray}\label{p2com3}
&&\sup _{||h||_{\dot{B}^0_{\infty,\infty}}\le 1}
\int_{\R^n}\sum_j\sum_{|t-j|\le 3}[{\cal{R}}\Delta^{1/4}(\Delta^{1/4}Q_j {\cal{R}} u^{j-4})-
\nabla(Q_j{\cal{R}}u^{j-4})]h_t dx\\
&&\simeq \sup _{||h||_{\dot{B}^0_{\infty,\infty}}\le 1}\int_{\R^n}\sum_j \sum_{|t-j|\le 3} {\cal{R}}u^{j-4}
[{\cal{R}}\Delta^{1/4} h_t \Delta^{1/4} Q_j)- \nabla h_t ) Q_j ]dx\nonumber\\
&&\simeq \sup _{||h||_{\dot{B}^0_{\infty,\infty}}\le 1}\int_{\R^n}\sum_j \sum_{|t-j|\le 3}
{\cal{F}}[{\cal{R}}u^{j-4}](\xi)(\int_{\R^n}{\cal{F}}[Q_j](\zeta){\cal{F}}[{\cal{R}}\Delta^{1/4} h_t](\xi-\zeta)\nonumber\\
&&~~~~ 
\left(|\zeta|^{1/2}- |\xi-\zeta|^{1/2}\right)d\zeta)d\xi\,.\nonumber
\end{eqnarray}
Now  we can proceed exactly as in \rec{p2com2} and get
 \begin{eqnarray*}
&&\sup _{||h||_{\dot{B}^0_{\infty,\infty}}\le 1}\int_{\R^n}\sum_j\sum_{|t-j|\le 3} [{\cal{R}}\Delta^{1/4}(\Delta^{1/4}Q_j {\cal{R}} u^{j-4})-
\nabla(Q_j{\cal{R}}u^{j-4})]h_t dx\\
&&\le C||Q||_{\dot{H}^{1/2}}||u||_{\dot{H}^{1/2}}\,.\hfill\Box
\end{eqnarray*}
\par\medskip
From Theorem \ref{comm3} and the duality between ${\cal{H}}^1$ and $BMO$ we obtain Theorem \ref{commestintr2}\,.
  
{\bf Proof  of Theorem \ref{commestintr2}.} It follows from Theorem \ref{comm3} and the duality between ${\cal{H}}^1$ and $BMO$\,.\hfill$\Box$
\begin{Lemma}\label{lemmareg} 
Let $u\in \dot H^{1/2}(\R^n)$, then ${\cal{R}}(\Delta^{1/4} u \cdot {\cal{R}}\Delta^{1/4} u)\in {\cal{H}}^1\,$ and
$$
||{\cal{R}}(\Delta^{1/4} u \cdot {\cal{R}}\Delta^{1/4} u)||_{{\cal{H}}^1}\le C||u||_{\dot{H}^{1/2}}^2\,.$$
\end{Lemma}
{\bf Proof of lemma~\ref{lemmareg}.} Since ${\cal{R}}\colon{\cal{H}}^1\to {\cal{H}}^1$, 
it  is enough to verify that  $\Delta^{1/4} u \cdot {\cal{R}}\Delta^{1/4} u\in {\cal{H}}^1\,.$\par
$\bullet $ Estimate of $\Pi_1(\Delta^{1/4} u\cdot {\cal{R}}\Delta^{1/4} u)$
\begin{eqnarray}\label{pi1n}
&&||\Pi_1(\Delta^{1/4} u \cdot {\cal{R}}\Delta^{1/4} u)||_{{\cal{H}}^1}=\int_{\R^n}(\sum_{j=-\infty}^{+\infty}[\Delta^{1/4} u_j({\cal{R}}\Delta^{1/4} u)^{j-4}]^2)^{1/2} dx\nonumber\\
&&\le
\int_{\R^n}\sup_j|({\cal{R}}\Delta^{1/4} u)^{j-4})|(\sum_{j=0}^{+\infty}[\Delta^{1/4} u_j]^2)^{1/2} dx\\
&&\le (\int_{\R^n} |M({\cal{R}}\Delta^{1/4} u)|^2dx)^{1/2}(\int_{\R^n}\sum_{j=-\infty}^{+\infty}[\Delta^{1/4} u_j]^2 dx)^{1/2}\nonumber\\
&&\le C||u||^2_{\dot{H}^{1/2}}\nonumber
\end{eqnarray}
\par
The estimate of the ${\cal{H}}^1$ norm of $\Pi_2(\Delta^{1/4} u \cdot {\cal{R}}\Delta^{1/4} u)$ is similar to \rec{pi1n}\,.

$\bullet $ Estimate of $\Pi_3(\Delta^{1/4} u \cdot {\cal{R}}\Delta^{1/4} u)$\,.
 \begin{eqnarray}\label{pi3n}
&&\sup _{||h||_{\dot{B}^0_{\infty,\infty}}\le 1}\int_{\R^n}\sum_{j}\sum_{|k-j|\le 3}\Delta^{1/4} u_j {\cal{R}}(\Delta^{1/4} u_k)[h^{j-6}+\sum_{t-=j-5}^{j+6}h_t]dx\nonumber\\
&&=\sup _{||h||_{\dot{B}^0_{\infty,\infty}}\le 1}\int_{\R^n}\sum_{j}\sum_{|k-j|\le 3}\left[\Delta^{1/4} u_j {\cal{R}}(\Delta^{1/4} u_k)-u_j\nabla u_k+\frac{1}{2}\nabla(u_j u_k)\right]\nonumber\\
&& ~~~~~~[h^{j-6}+\sum_{t-=j-5}^{j+6}h_t]dx
\end{eqnarray}
We only estimate the terms with $h^{j-6}$, being the estimates with $h_t$ similar \,.
\par
\begin{eqnarray*}
&& \sup _{||h||_{\dot{B}^0_{\infty,\infty}}\le 1}\int_{\R^n}\sum_{j}\sum_{|k-j|\le 3}(\Delta^{1/4} u_j {\cal{R}}(\Delta^{1/4} u_k)-u_j\nabla u_k)h^{j-6}dx\\
&&  \sup _{||h||_{\dot{B}^0_{\infty,\infty}}\le 1}\int_{\R^n}\sum_{j}\sum_{|k-j|\le 3}{\cal {F}}[h^{j-6}](x)\left(\int_{\R^n} {\cal {F}}[u_j]{\cal {F}}[{\cal{R}}\Delta^{1/4} u_k][|y|^{1/2}-|x-y|^{1/2}] dy\right)dx\\
&& ~\mbox{by arguing as in \rec{p2com2}}\\
&&\le C||u||_{\dot{H}^{1/2}}^2
\end{eqnarray*}
Finally we also have
\begin{eqnarray*}
&& \sup _{||h||_{\dot{B}^0_{\infty,\infty}}\le 1}\int_{\R^n}\sum_{j}\sum_{|k-j|\le 3}\frac{1}{2}\nabla(u_j u_k)h^{j-6} dx \\
&&~~~=\sup _{||h||_{\dot{B}^0_{\infty,\infty}}\le 1}\int_{\R^n}\sum_{j}\sum_{|k-j|\le 3}\frac{1}{2}(u_j u_k)\nabla h^{j-6} dx\\
&&~~~\le C\sup _{||h||_{\dot{B}^0_{\infty,\infty}}\le 1}||h||_{\dot{B}^0_{\infty,\infty}}
 \int_{\R^n}\sum_{j}\sum_{|k-j|\le 3}  2^{j}u_j u_k dx\\
&&~~~\le C ( \int_{\R^n}\sum_{j} 2^ju_j^2 dx)^{1/2}= C ||u||_{\dot{H}^{1/2}}^2\,.~~~~~~~~~\hfill\mbox{$\Box$}
\end{eqnarray*}
 
 We get the following result
 \begin{Corollary}
 Let $n\in \dot H^{1/2}(\R^n,S^{m-1})$. Then
 \begin{equation}\label{nest}
 \Delta^{1/4}[n\cdot \Delta^{1/4} n] \in {\cal{H}}^{1}\,.
 \end{equation}
 \end{Corollary}
 {\bf Proof.} 
 Since $n\cdot \nabla n=0$  we can write
 \begin{eqnarray*}
  \Delta^{1/4}[n\cdot \Delta^{1/4} n]&= &\Delta^{1/4}[n\cdot \Delta^{1/4} n]-{\cal{R}}(n\cdot\nabla n)
  + {\cal{R}}[\Delta^{1/4}n \cdot {\cal{R}}\Delta^{1/4} n]\\
  &&~~~~ -{\cal{R}}[\Delta^{1/4}n\cdot {\cal{R}}\Delta^{1/4} n]\\
  &=& S(n\cdot,n)--{\cal{R}}[\Delta^{1/4}n\cdot {\cal{R}}\Delta^{1/4} n]\,.
  \end{eqnarray*}
  The estimate \rec{nest} is a consequence of Theorem \ref{comm3}  and
    Lemma \ref{lemmareg}, which respectively imply that $S(n\cdot,n)\in {\cal{H}}^{1}$ and ${\cal{R}}(\Delta^{1/4}n\cdot {\cal{R}}\Delta^{1/4} n)\in {\cal{H}}^{1}\,.$~\hfill$\Box$
    
    
  \section{Geometric localization properties of the $\dot{H}^{1/2}-$norm on the real line.}

  In the next Theorem we  show that the $\dot{H}^{1.2}([a,b])$ norm ($-\infty\le a<b\le+\infty$) can be localized in space.
  This result, besides being of independent interest, will be used in Section \ref{regulsection} for suitable localization estimates.
    For simplicity we will suppose that $[a,b]=[-1,1]$.
\begin{Theorem}\label{localization}{\bf [Localization of $H^{1/2}((-1,1))$ norm]}
Let $u\in H^{1/2}((-1,1))\,.$  Then for some $C>0$ we have
$$||u||^2_{\dot{H}^{1/2}((-1,1))}\simeq\sum_{j=-\infty}^{0}||u||^2_{\dot{H}^{1/2}(A_j)}
$$
where $A_j=B_{2^{j+1}}\setminus B_{2^{j-1}}\,.$ 
\end{Theorem}
{\noindent\bf Proof.}  We set for every $i\in \Z$, $A^\prime_i=B_{2^{i}}\setminus B_{2^{i-1}}$  and
$\bar{u}^\prime_i=|A_i^\prime|^{-1}\int_{A^\prime_i} u(x)\,dx$ (i.e. the mean value of $u$ on the annulus
  $A_i^\prime$).
We have
\begin{eqnarray}\label{decomp}
&& ||u||^2_{\dot{H}^{1/2}((-1,1))} \simeq\int_{[-1,1]}\int_{[-1,1]}\frac{|u(x)-u(y)|^2}{|x-y|^2}dx dy\\
&&=
\sum_{i,j=-\infty}^0\int_{A^\prime_i}\int_{A^{\prime}_j}\frac{|u(x)-u(y)|^2}{|x-y|^2}dx dy\nonumber \\
&&=
\sum_{i=-\infty}^0\int_{A^\prime_i}\int_{A^{\prime}_i}\frac{|u(x)-u(y)|^2}{|x-y|^2}dx dy\nonumber\\
&&+ 2\sum_{j=-\infty}^0\sum_{i>j+1}\int_{A^\prime_i}\int_{A^{\prime}_j}\frac{|u(x)-u(y)|^2}{|x-y|^2}dx dy\nonumber\\&&
+2\sum_{j=-\infty}^0\int_{A^\prime_j}\int_{A^{\prime}_{j+1}}\frac{|u(x)-u(y)|^2}{|x-y|^2}dx dy\,.\nonumber
\end{eqnarray}
We first observe that
\begin{equation}\label{m3}
\sum_{i,j=-\infty}^0\int_{A^\prime_i}\int_{A^{\prime}_j}\frac{|u(x)-u(y)|^2}{|x-y|^2}dx dy\le
\sum_{i,j=-\infty}^0\int_{A_i}\int_{A_j}\frac{|u(x)-u(y)|^2}{|x-y|^2}dx dy\end{equation}
and
\begin{equation}\label{m4}
\sum_{j=-\infty}^0\int_{A^\prime_j}\int_{A^{\prime}_{j+1}}\frac{|u(x)-u(y)|^2}{|x-y|^2}dx dy\le
\sum_{j=-\infty}^0\int_{A_j}\int_{A_{j}}\frac{|u(x)-u(y)|^2}{|x-y|^2}dx dy\,.\end{equation}
It remains to estimate the term $\sum_{j=-\infty}^0\sum_{i>j+1}\int_{A^\prime_i}\int_{A^{\prime}_j}\frac{|u(x)-u(y)|^2}{|x-y|^2}dx dy$ 
 in (\ref{decomp}).\par  We have
\begin{eqnarray*}
&&
\sum_{j=-\infty}^0\sum_{i>j+1}\int_{A^\prime_i}\int_{A^{\prime}_j}\frac{|u(x)-u(y)|^2}{|x-y|^2}dx dy\\
&& \le C \sum_{j=-\infty}^0\sum_{i\ge j+2}2^{-2i}\int_{A^\prime_i}\int_{A^{\prime}_j}{|u(x)-u(y)|^2} dx dy\\
&& \le 
C (\sum_{j=-\infty}^0\sum_{i\ge j+2}2^{-2i}\int_{A^\prime_i}\int_{A^{\prime}_j}|\bar{u}^\prime_i-\bar{u}^\prime_j |^2dx dy\\
&& +
 \sum_{j=-\infty}^0\sum_{i\ge j+2}2^{-2i}\int_{A^\prime_i}\int_{A^{\prime}_j}{|u(x)-\bar{u}^\prime_i|^2} dx dy\\
 && +
  \sum_{j=-\infty}^0\sum_{i\ge j+2}2^{-2i}\int_{A^\prime_i}\int_{A^{\prime}_j}{|u(y)-\bar{u}^\prime_j|^2} dx dy)\\
  &&
  \le C (\sum_{j=-\infty}^0\sum_{i\ge j+2}2^{-2i} 2^{i+j}|\bar{u}^\prime_i-\bar{u}^\prime_j |^2\\
&& +
 \sum_{j=-\infty}^0\sum_{i\ge j+2}2^{-2i}2^{j}\int_{A^\prime_i} {|u(x)-\bar{u}^\prime_i|^2} dx \\
 && +
  \sum_{j=-\infty}^0\sum_{i\ge j+2}2^{-2i} 2^{i}\int_{A^{\prime}_j}{|u(y)-\bar u^\prime_j|^2} dy)\,.
  \end{eqnarray*}
  
 $\bullet$ Estimate of $\sum_{j=-\infty}^0\sum_{i\ge j+2}2^{-2i}2^{j}\int_{A^\prime_i} {|u(x)-\bar{u}^\prime_i|^2} dx$\,.
 
  \begin{eqnarray}\label{m1}
  && \sum_{j=-\infty}^0\sum_{i\ge j+2}2^{-2i}2^{j}\int_{A^\prime_i} {|u(x)-\bar{u}^\prime_i|^2} dx\\
  &&=\sum_{i=-\infty}^0\sum_{j\le i-2}2^{-2i}2^{j}\int_{A^\prime_i} {|u(x)-\bar{u}^\prime_i|^2} dx\nonumber\\
  &&=\sum_{i=-\infty}^0 2^{-2i}\int_{A^\prime_i} {|u(x)-\bar u^\prime_i|^2} dx(\sum_{j\le i-2}2^j)\nonumber\\
  &&\le C\sum_{i=-\infty}^0 |A^\prime_i|^{-1}\int_{A^\prime_i} {|u(x)- \bar u^\prime_i|^2} dx\nonumber\\
  &&\le C \sum_{i=-\infty}^0\int_{A^\prime_i}\int_{A^\prime_i} \frac {|u(x)-u(y)|^2} {|x-y|^2 }dx dy \,.\nonumber
   \end{eqnarray}
   In the last inequality we use the fact that for every $i$ it holds
      \begin{eqnarray*}
   &&
   |A^\prime_i|^{-1}\int_{A^\prime_i} {|u(x)- \bar u^\prime_i|^2} dx\\
   &&\le  |A^\prime_i|^{-1}\int_{A^\prime_i} {|u(x)- |A^\prime_i|^{-1}\int_{A^\prime_i}u(y) dy|^2} dx\\
   && \le |A^\prime_i|^{-2}\int_{A^\prime_i}\int_{A^\prime_i}|u(x)-u(y)|^2 dx dy\\
   && \le C\int_{A^\prime_i}\int_{A^\prime_i}\frac{|u(x)-u(y)|^2 }{|x-y|^2 }dyx dy\,.
   \end{eqnarray*}
  $\bullet$   Estimate of $\sum_{j=-\infty}^0\sum_{i\ge j+2}2^{-2i}2^{j}\int_{A^\prime_i} {|u(y)-\bar{u}^\prime_j|^2} du$
    \begin{eqnarray}\label{m2}
  && \sum_{j=-\infty}^0\sum_{i\ge j+2}2^{-i}\int_{A^\prime_j} {|u(y)-\bar u^\prime_j|^2} dy\\
     &&=\sum_{j=-\infty}^0 \int_{A^\prime_j} {|u(x)-\bar u_j|^2} dx(\sum_{i\le j+2}2^{-i})\nonumber\\
  &&=\frac{1}{2}\sum_{j=-\infty}^0 2^{-j}\int_{A^\prime_j} {|u(x)-\bar{u}^\prime_j|^2} dy\nonumber\\
  &&\le C \sum_{j=-\infty}^0  \int_{A^\prime_j}\int_{A^\prime_j} \frac {|u(x)-u(y)|^2} {|x-y|^2} dx dy \,.\nonumber
   \end{eqnarray}
   $\bullet$ Estimate of $\sum_{j=-\infty}^0\sum_{i\ge j+2}2^{-2i} 2^{i+j}|\bar{u}^\prime_i-\bar{u}^\prime_j |^2\,.$
  We first  observe that
   $$
   |\bar{u}^\prime_i-\bar{ u}^\prime_j|^2\le (i-j)\sum_{j}^{i-1}|\bar u^\prime_{\ell+1}-\bar u^\prime_\ell|^2
   $$
   and
   $$
   |\bar u_{\ell+1}-\bar u_\ell|^2\le |A_\ell|^{-1}\int_{A_\ell}|u-\bar {u}_\ell|^2 dx,
   $$
   where
   $\bar {u}_\ell=|A_\ell|^{-1}\int_{A_\ell}u(x)\,dx\,.$
   \par 
   We set
   $a_\ell=|A_\ell|^{-1}\int_{A_\ell}|u-\bar {u}_\ell|^2 dx$\,.
   We have
   \begin{eqnarray*}
 &&  \sum_{j=-\infty}^0\sum_{i\ge j+2}2^{-2i} 2^{i+j}|\bar{u}^\prime_i-\bar{u}^\prime_j |^2\\
   &&\le  \sum_{j=-\infty}^0\sum_{i\ge j+2} (i-j)2^{j-i}\sum_j^{i-1} a_\ell \leq 
   \sum_{\ell=-\infty}^0 a_\ell \sum_{j=-\infty}^{\ell} \sum_{i-j\ge \ell+1-j}  (i-j)2^{j-i}\,.
   \end{eqnarray*}
   We observe that
   \begin{eqnarray} 
   &&\sum_{i-j\ge \ell+1-j}  (i-j)2^{j-i}\le \int_{\ell+1-j }^{+\infty}2^{-x} x dx =2^{-(\ell+1-j)}(\ell+2-j)\,
   \end{eqnarray}
   and
 $$
       \sum_{j=-\infty}^{\ell}2^{-(\ell+1-j)}(\ell+2-j)\le \int_1^{+\infty}2^{-t}(t+1)\,dx\le C,
   $$
   for some constant $C$ independent on $\ell\,.$
 Therefore  we get
    \begin{eqnarray}\label{m6}
   && \sum_{j=-\infty}^0\sum_{i\ge j+2}2^{-2i} 2^{i+j}|\bar{u}^\prime_i-\bar{u}^\prime_j |^2\\
    &&\le   \sum_{j=-\infty}^0\sum_{i\ge j+2} (i-j)2^{j-i}\sum_j^{i-1} a_\ell\le  C \sum_{\ell=-\infty}^0 a_\ell\le C \sum_{\ell=-\infty}^0\int_{A_\ell}\int_{A_\ell}\frac{|u(x)-u(y)|^2}{|x-y|^2}\,dx dy\,.\nonumber
    \end{eqnarray}
   By combining \rec{m3},\rec{m4},\rec{m1},\rec{m2} and \rec{m6} we finally obtain
   $$
   ||u||^2_{\dot{H}^{1/2}((-1,1))}\lesssim   \sum_{\ell=-\infty}^0 ||u||^2_{\dot{H}^{1/2}(A_\ell)} \,.
   $$
   Next we show that
   \begin{equation}\label{other}
    \sum_{\ell=-\infty}^0 ||u||^2_{\dot{H}^{1/2}(A_\ell)}\lesssim  ||u||^2_{\dot{H}^{1/2}((-1,1))}\,.
    \end{equation}
    We observe that
    for every $\ell$ we have $A_\ell=C_\ell\cup D_\ell$ where $C_\ell =B_{2^{\ell+1}}\setminus B_{2^{\ell}}$ and
    $ D_\ell=B_{2^{\ell}}\setminus B_{2^{\ell-1}}\,.$ Thus
    \begin{eqnarray*}
&& ||u||^2_{\dot{H}^{1/2}(A_\ell)}=\int_{C_{\ell}}\int_{C_\ell} \frac{|u(x)-u(y)|^2}{|x-y|^2}dx dy\\
&&~~+\int_{D_{\ell,h}}\int_{D_\ell} \frac{|u(x)-u(y)|^2}{|x-y|^2}dx dy+2\int_{D_{\ell,h}}\int_{C_\ell} \frac{|u(x)-u(y)|^2}{|x-y|^2}dx dy\,.
\end{eqnarray*}
Since   $\cup_{\ell} (C_\ell\times C_\ell)$, $\cup_{\ell} (D_\ell\times C_\ell)$ and $\cup_{\ell} (D_\ell\times C_\ell)$ are disjoint unions contained in $[0,1]\times [0,1]$ we
have
$$
\sum_\ell\int_{C_{\ell}}\int_{C_\ell} \frac{|u(x)-u(y)|^2}{|x-y|^2}dx dy\le \int_{[-1,1]}\int_{[-1,1]}\frac{|u(x)-u(y)|^2}{|x-y|^2}dx dy\,;$$
$$
\sum_\ell\int_{D_{\ell,h}}\int_{C_\ell} \frac{|u(x)-u(y)|^2}{|x-y|^2}dx dy\le \int_{[-1,1]}\int_{[-1,1]}\frac{|u(x)-u(y)|^2}{|x-y|^2}dx dy\,;$$
$$
\sum_\ell\int_{D_{\ell,h}}\int_{D_\ell} \frac{|u(x)-u(y)|^2}{|x-y|^2}dx dy\le \int_{[-1,1]}\int_{[-1,1]}\frac{|u(x)-u(y)|^2}{|x-y|^2}dx dy\,.$$
It follows that
$$ \sum_{\ell=-\infty}^0 ||u||^2_{\dot{H}^{1/2}(A_\ell)}\le \bar C \int_{[-1,1]}\int_{[-1,1]}\frac{|u(x)-u(y)|^2}{|x-y|^2}dx dy= \bar C ||u||^2_{\dot{H}^{1/2}((-1,1))}\,$$
and we can conclude.~\hfill$\Box$
\begin{Remark}{\rm By analogous computations one can show that for all $r>0$ we have

$$||u||^2_{\dot{H}^{1/2}(\R)} \simeq\sum_{j=-\infty}^{+\infty}||u||^2_{\dot{H}^{1/2}(A_j^r)}
$$
where $A^r_j=B_{2^{j+1}r}\setminus B_{2^{j-1}r}\,,$ where the equivalence constants do not depend on $r$.} 

\end{Remark}
   \par\bigskip
   Next we compare the $\dot{H}^{1/2}$ norm of $\Delta^{-1/4}( M\Delta^{1/4}u)$ with the $L^{2}$ norm of   $M\Delta^{1/4}u$,  where $u\in\dot{H}^{1/2}(\R)$ and   
   $M\in \dot{H}^{1/2}(\R,{\cal{M}}_{p\times m}(\R))$, $p\ge 0\,.$\par

      \begin{Lemma}\label{rotation}
   Let  $M\in \dot{H}^{1/2}(\R,{\cal{M}}_{p\times m}(\R))$, $m\ge1, p\ge 1$, and $u\in \dot{H}^{1/2}(\R)$.
    Then   there exist $C_1>0$, $C_2>0$ and $n_0\in{\N}$, independent of $u$ and $M$, such that, for any $r\in (0,1)$, $n>n_0$
    and any $x_0\in{\R}$, we have 
     
   \begin{eqnarray*}
   ||\Delta^{-1/4}( M\Delta^{1/4}u)||^2_{\dot{H}^{1/2}(B_r(x_0))}
   &&
   \ge C_1\int_{B_{r/2^n}(x_0)}|M\Delta^{1/4} u|^2 dx\\
   &&
    - C_2\sum_{h=-n}^{+\infty}2^{-h}\int_{B_{2^{h}r}(x_0)\setminus B_{2^{h-1}r}(x_0)}|M\Delta^{1/4} u|^2\ dx\,.
   \end{eqnarray*}
   \end{Lemma}

  {\bf Proof of lemma~\ref{rotation}.}
  We write
  \begin{eqnarray*}
  &&
  \Delta^{-1/4}( M\Delta^{1/4} u)= \Delta^{-1/4}(\11_{|x|\le r/2^n} M\Delta^{1/4} u)+
  \Delta^{-1/4}((1-\11_{|x|\le r/2^n}) M\Delta^{1/4} u)\,,
  \end{eqnarray*}
  where $n>0$ is large enough(the threshold will be determined later in the proof)\,.\par
  For any $\rho\ge 0$, we denote by
  $\11_{|x|\le\rho}$ and $\11_{ \rho\le |x|}$
 the characteristic functions of the sets of point $x\in{\R}$ respectively where $|x|\le\rho $ and $|x|\ge \rho$. For all $\rho\le\sigma$ we also denote by $\11_{\rho\le |x|\le \sigma}$ the characteristic function
 of the set $\{x\in{\R}\ ;\ \rho\le |x|\le \sigma\}$. We have
  \begin{equation}\label{estrot1}
  \begin{array}{l}
 \ds || \Delta^{-1/4}( M\Delta^{1/4} u)||_{\dot{H}^{1/2}(B_r)}\ge
 || \Delta^{-1/4}(\11_{r/2^n} M\Delta^{1/4} u)||_{\dot{H}^{1/2}(B_r)}\\[5mm]
 \ds~~~~~~~~~~~~~~~-|| \Delta^{-1/4}((1-\11_{|x|\le r/2^n}) M\Delta^{1/4} u)||_{\dot{H}^{1/2}(B_r)}\\[5mm]
 \ds \ge
  || \Delta^{-1/4}(\11_{r/2^n} M\Delta^{1/4} u)||_{\dot{H}^{1/2}(B_r)}-
  || \Delta^{-1/4}(\11_{r/2^n\le |x|\le 4r} M\Delta^{1/4} u)||_{\dot{H}^{1/2}(B_r)}\\[5mm]
 \ds ~~~~~~~~~~~~~~-
  || \Delta^{-1/4}(\11_{|x|\ge 4r} M\Delta^{1/4} u)||_{\dot{H}^{1/2}(B_r)}\\[5mm]
  \ds\ge
  || \Delta^{-1/4}(\11_{|x|\le r/2^n} M\Delta^{1/4} u)||_{\dot{H}^{1/2}(B_r)}-
  || \Delta^{-1/4}(\11_{r/2^n\le |x|\le 4r} M\Delta^{1/4} u)||_{\dot{H}^{1/2}(\R)} \\[5mm]
  ~~~~~~~~~~~~~~-
  || \Delta^{-1/4}(\11_{|x|\ge 4r} M\Delta^{1/4} u)||_{\dot{H}^{1/2}(B_r)}\quad.
  \end{array}
  \end{equation}\par
  
  We estimate of the last three terms  in \rec{estrot1}\,.
  
  \medskip
  
  $\bullet$  Estimate of $|| \Delta^{-1/4}(\11_{r/2^n\le |x|\le 4r} M\Delta^{1/4} u)||_{\dot{H}^{1/2}(\R)} $.\par
  \begin{equation}\label{estrot2}
  \begin{array}{l}
  
 \ds || \Delta^{-1/4}(\11_{r/2^n\le |x|\le 4r} M\Delta^{1/4} u)||^2_{\dot{H}^{1/2}(\R)} =
  \int_{r/2^n\le |x|\le 4r} | M\Delta^{1/4} u|^2\ dx\\
  
  \ds\quad=\sum_{h=-n}^1\int_{2^{h}r\le |x|\le 2^{h+1} r}| M\Delta^{1/4} u|^2\ dx\,.
  \end{array}
  \end{equation}
  
   $\bullet$ Estimate of $ || \Delta^{-1/4}(\11_{|\xi|\ge 4r} M\Delta^{1/4} u)||_{\dot{H}^{1/2}(B_r)}\,.$ We set
   $$g:=\11_{|x|\ge 4r} M\Delta^{1/4} u\quad.$$
   With this notation we have
   \begin{equation}
   \label{estrotz1}
   \begin{array}{l}
  \ds  || \Delta^{-1/4}(\11_{|x|\ge 4r} M\Delta^{1/4} u)||_{\dot{H}^{1/2}(B_r)}^2\\[5mm]
    =
  \ds  \int_{B_r}\int_{B_r}\frac{|(\frac{1}{|x|^2}\star g)(t)-(\frac{1}{|x|^2}\star g)(s)|^2}{|x-y|^2} dt\ ds\\[5mm]
   \ds =
     \int_{B_r}\int_{B_r}\frac{1}{|t-s|^2}\left(\int_{|x|\ge 4r} g(x)(\frac{1}{|t-x|^{1/2}}-\frac{1}{|s-x|^{1/2}})dx \right)^2 dt\ ds\\[5mm]
    ~~\mbox{by Mean Value Theorem}\\[5mm]
    
\ds     \le C \int_{B_r}\int_{B_r}\left(\int_{|x|\ge 4r} |g(x)|\max(\frac{1}{|t-x|^{3/2}},\frac{1}{|s-x|^{3/2}})dx\right)^2 dt\ ds\\[5mm]
 
\ds  \le C \int_{B_r}\int_{B_r}\left(\sum_{h=4}^{+\infty}\int_{2^{h}r\le |x|\le 2^{h+1}r} |g(x)|\max(\frac{1}{|t-x|^{3/2}},\frac{1}{|s-x|^{3/2}})dx\right)^2 dt\ ds\\[5mm]
   
\ds   \le C \int_{B_r}\int_{B_r}\left(\sum_{h=4}^{+\infty}\int_{2^{h}r\le |x|\le 2^{h+1}r} |g(x)| 2^{-3/2h}r^{-3/2}d\xi\right)^2 dt\ ds\\[5mm]
  ~~\mbox{by H\"older Inequality}\\[5mm]
   
 \ds  \le C \int_{B_r}\int_{B_r}\left(\sum_{h=4}^{+\infty}2^{-h}r^{-1}(\int_{2^{h+1}r\le |x|\le 2^{h+1}r} |g(x)| ^2\ dx)^{1/2}\right)^2 dt\ ds\\[5mm]
   ~~\mbox{by Cauchy-Schwartz Inequality}\\[5mm]
   
\ds   \le 
   C(\sum_{h=4}^{+\infty} 2^{-h})\left(\sum_{h=4}^{+\infty} 2^{-h}\int_{B_{2^{h+1}r}(x_0)\setminus B_{2^{h}r}(x_0)}|M\Delta^{1/4} u|^2\ dx\right)\\[5mm]
   \ds \le C\left(\sum_{h=4}^{+\infty} 2^{-h}\int_{B_{2^{h+1}r}(x_0)\setminus B_{2^{h}r}(x_0)}|M\Delta^{1/4} u|^2\ dx\right)\,.
   \end{array}
   \end{equation}
   
   \medskip
   
     $\bullet$  Estimate of  $|| \Delta^{-1/4}(\11_{|x|\le r/2^n} M\Delta^{1/4} u)||_{\dot{H}^{1/2}(B_r)}$\,.\par
     We set
     $$A_h^r:=\{x~:~2^{h-1}r\le|x|\le 2^{h+1} r\}\quad .$$
     By  Localization Theorem \ref{localization} there exists a constant $\tilde C>0$ (independent  on $r$ ) such that
     \begin{equation}
     \label{estrotz2}
     \begin{array}{rl}
    \ds || \Delta^{-1/4}(\11_{|x|\le r/2^n} M\Delta^{1/4} u)||^2_{\dot{H}^{1/2}(\R)} 
   &\ds\le \tilde C\sum_{h=-\infty}^{+\infty} || \Delta^{-1/4}(\11_{|x|\le r/2^n} M\Delta^{1/4} u)||^2_{\dot{H}^{1/2}(A^r_h)}\\[5mm]
   &\ds\le 
    \tilde C\   || \Delta^{-1/4}(\11_{|x|\le r/2^n} M\Delta^{1/4} u)||^2_{\dot{H}^{1/2}(B_r)}\\[5mm]
   &\ds+
   \tilde C
   \sum_{h=0}^{+\infty} || \Delta^{-1/4}(\11_{|x|\le r/2^n} M\Delta^{1/4} u)||^2_{\dot{H}^{1/2}(A^r_h)}\,.
   \end{array}
   \end{equation}
   $\bullet$ Estimate of 
   $\sum_{h=0}^{+\infty} || \Delta^{-1/4}(\11_{|x|\le r/2^n} M\Delta^{1/4} u)||^2_{\dot{H}^{1/2}(A^r_h)}$.
   We set now
   \[
   f(x):=\11_{|x|\le r/2^n}\ (M\Delta^{1/4}u)\,.
   \]
   Using this notation we have
   \begin{equation}
   \label{estrotz3}
   \begin{array}{l}
   \ds\sum_{h=0}^{+\infty} || \Delta^{-1/4}(\11_{|x|\le r/2^n} M\Delta^{1/4} u)||^2_{\dot{H}^{1/2}(A^r_h)}\\[5mm]
   \ds\le 
 \sum_{h=0}^{+\infty} \int_{A_h^r}\int_{A_h^r}\left(\int_{|x|\le r/2^n} |f(x)|(|\frac{1}{|t-x|^{1/2}}-\frac{1}{|s-x|^{1/2}}|)dx\right)^2 dt\ ds\\[5mm]
     ~~\mbox{by Mean Value Theorem}\\[5mm]
     \ds\le C \sum_{h=0}^{+\infty} \int_{A_h^r}\int_{A_h^r}\left(\int_{|x|\le r/2^n} |f(x)|\max(\frac{1}{|t-x|^{3/2}},\frac{1}{|s-x|^{3/2}})d\xi\right)^2 dt\ ds\\[5mm]
\ds
\le C\sum_{h=0}^{+\infty} \int_{A_h^r}\int_{A_h^r} \max(\frac{1}{|t|^3},\frac{1}{|s|^3})\frac{r}{2^n}
(\int_{|x|\le r/2^n}|f(x)|^2\ dx) dt\ ds\\[5mm]
\ds=
\frac{C}{2^{n}}\sum_{h=0}^{+\infty} 2^{-h}(\int_{|x|\le r/2^n}|f(x)|^2\ dx)\le\frac{C}{2^n}\int_{B_{r/2^n}(x_0)}|M\Delta^{1/4}u|^2\ dx
\end{array}
\end{equation}
If  $n$ is large enough in such a way that $C\,{\tilde C}/2^n<1/2$, we get, combining (\ref{estrot1}), (\ref{estrot2}), (\ref{estrotz1}), (\ref{estrotz2}) and (\ref{estrotz3}), for some $C_1,C_2$ positive,
\begin{equation}\label{estrot3}
\begin{array}{l}
||\Delta^{-1/4}(M\Delta^{1/4}u)||^2_{\dot{H}^{1/2}(B_r)}\\[5mm]
\displaystyle\quad\quad \ge C_1 \int_{B_{r/2^n}}| M\Delta^{1/4} u|^2\ dx\\[5mm]
\ds
  \displaystyle \quad\quad  - C_2\sum_{h=-n}^{+\infty}2^{-h}\int_{B_{2^{h+1}r}\setminus B_{2^{h}r}}|M\Delta^{1/4} u|^2\ dx \,,
\end{array}
   \end{equation}
   which ends the proof of the lemma.\hfill $\Box$
   \par\medskip
   In the following Lemma we compare the $\dot H^{1/2}$ norm of $w= \Delta^{-1/4}( M\Delta^{1/4} u)$ in the annuli
  $A_h=B_{2^{h+1}}(x_0)\setminus B_{2^{h-1}}(x_0)$ with the $L^2$ norm in the same annuli of $M\Delta^{1/4} u$.
  Such a result will be used in the following Section for suitable
 {\em localization} estimates.\par\medskip
\begin{Lemma}\label{rotation2}
Let  $M\in \dot{H}^{1/2}(\R,{\cal{M}}_{p\times m}(\R))$, $m\ge1, p\ge1$, and $u\in \dot H^{1/2}(\R)$. Then there exists $C>0$ such that for every $\gamma\in (0,1)$, for all $n\ge  n_0\in{\N}$ ($n_0$ dependent on $\gamma$ and independent of $u$ and $M$), for every $k\in\Z $,  and 
   any $x_0\in{\R}$,    we have
\begin{eqnarray}\label{anelli}
&&\sum_{h=k}^{+\infty} 2^{k-h}||\Delta^{-1/4}( M\Delta^{1/4} u)||^2_{\dot{H}^{1/2}(B_{2^{h+1}}(x_0)\setminus B_{2^{h-1}}(x_0))}\le 
\gamma \int_{B_{2^{k-n}}(x_0)}| M\Delta^{1/4} u|^2 d\xi\nonumber\\
&&
+ \sum_{h=k-n}^{+\infty}2^{\frac{k-h}{2}}\int_{B_{2^{h+1}}(x_0)\setminus B_{2^{h-1}}(x_0)}| M\Delta^{1/4} u|^2 d\xi\,.\nonumber
\end{eqnarray}
\end{Lemma}
{\bf Proof.}
Given $h \in\Z$ and $\ell\ge 3$ we set $A_h=B_{2^{h+1}}(x_0)\setminus B_{2^{h-1}}(x_0)$ and $D_{\ell,h}=B_{2^{h+\ell}}(x_0)\setminus B_{2^{h-\ell}}(x_0)$\,.
For simplicity of notations we suppose that $x_0=0$ but all the following estimates will be independent on $x_0$.
   \par
   We fix $\gamma\in (0,1)\,.$\par
We have
\begin{eqnarray*}
&&
||w||^2_{\dot{H}^{1/2}(A_h)}=\int_{A_h}\int_{A_h}\frac{|w(x)-w(y)|^2}{|x-y|^2} dx dy\\ [5mm]
&& \le 2||\Delta^{-1/4}\11_{D_{\ell,h}} M\Delta^{1/4} u||^2_{\dot{H}^{1/2}(A_h)}+
2||\Delta^{-1/4}(1-\11_{D_{\ell,h}}) M\Delta^{1/4} u||^2_{\dot{H}^{1/2}(A_h)}\\[5mm]
&&\le 
2||\Delta^{-1/4}\11_{D_{\ell,h}} M\Delta^{1/4} u||^2_{\dot{H}^{1/2}(\R)} +
2||\Delta^{-1/4}(1-\11_{D_{\ell,h}}) M\Delta^{1/4} u||^2_{\dot{H}^{1/2}(A_h)}\,.
\end{eqnarray*}
 The constant $\ell$ will be determined later.\par\medskip
$\bullet$ Estimate of $||\Delta^{-1/4}\11_{D_{\ell,h}} M\Delta^{1/4} u||^2_{\dot{H}^{1/2}(\R)} \,.$
\begin{eqnarray}\label{estrotD}
&&
||\Delta^{-1/4}\11_{D_{\ell,h}} M\Delta^{1/4} u||^2_{\dot{H}^{1/2}(\R)} =\int_{D_{\ell,h}}|M\Delta^{1/4} u|^2dx\nonumber\\
&&=\sum_{s=h-\ell}^{h+\ell-1}\int_{B_{2^{s+1}}\setminus B_{2^{s}}}|M\Delta^{1/4} u|^2dx\,.
\end{eqnarray}
We multiply \rec{estrotD} by $2^{k-h}$ and we sum up from $h=k$ to $+\infty$ and get
\begin{eqnarray}\label{estrDbis}
&&
\sum_{h=k}^{+\infty} 2^{k-h}||\Delta^{-1/4}\11_{D_{\ell,h}} M\Delta^{1/4} u||^2_{\dot{H}^{1/2}(\R)}\le
C 2^{\ell}\sum_{h=k-\ell}^{+\infty}\int_{B_{2^{h+1}}\setminus B_{2^{h-1}}}|M\Delta^{1/4} u|^2dx\,.
\end{eqnarray}

$\bullet$ Estimate of $||\Delta^{-1/4}(1-\11_{D_{\ell,h}}) M\Delta^{1/4} u||^2_{\dot{H}^{1/2}(A_h)}$.\par
 We set
$g=(1-\11_{D_{\ell,h}}) M\Delta^{1/4} u$\,.
\begin{eqnarray}\label{estrot4}
&&||\Delta^{-1/4}(1-\11_{D_{\ell,h}}) M\Delta^{1/4} u||^2_{\dot{H}^{1/2}(A_h)}=
\int_{A_h}\int_{A_h}\frac{|(\frac{1}{|x|^2}\star g)(t)-(\frac{1}{|x|^2}\star g)(s)|^2}{|t-s|^2} dt ds\nonumber\\
&&\le
2\int_{A_h}\int_{A_h}\frac{1}{|t-s|^2}\left(\int_{|x|> 2^{\ell+h}}  g(x)(\frac{1}{|x-t|^{1/2}}-\frac{1}{|x-s|^{1/2}})dx\right)^2  dt ds \\\
&&
+2\int_{A_h}\int_{A_h}\frac{1}{|t-s|^2}\left(\int_{|x|<2^{h-\ell}}  g(x)(\frac{1}{|x-t|^{1/2}}-\frac{1}{|x-s|^{1/2}})dx\right)^2  dt ds\,.\nonumber\
\end{eqnarray}
We estimate the last two terms in \rec{estrot4}. \par

{\bf 1.} Estimate of $\int_{A_h}\int_{A_h}\frac{1}{|t-s|^2}\left(\int_{|x|> 2^{\ell+h}}  g(x)(\frac{1}{|x-t|^{1/2}}-\frac{1}{|x-s|^{1/2}})dx\right)^2  dt ds\,.$
\begin{eqnarray}\label{estrot4bbis}
&&
\int_{A_h}\int_{A_h}\frac{1}{|t-s|^2}\left(\int_{|x|> 2^{\ell+h}}  g(x)(\frac{1}{|x-t|^{1/2}}-\frac{1}{|x-s|^{1/2}})dx\right)^2  dt ds\nonumber\\
&&
\le C
\int_{A_h}\int_{A_h}\left(\sum_{s=h+\ell}^{\infty}\int_{2^s\le|x|\le 2^{s+1}}  g(x)\max(\frac{1}{|x-t|^{3/2}},\frac{1}{|x-s|^{3/2}})dx\right)^2  dt ds\nonumber
\\
&&\mbox{by H\"older Inequality}\\
&&
\le C \int_{A_h}\int_{A_h}\left(\sum_{s=h+\ell}^{\infty}2^{-s}\int_{2^s\le|x|\le 2^{s+1}} | g(x)|^2dx)^{1/2}\right)^2  dt ds\nonumber\\
&&~~~\mbox{by Cauchy-Schwartz Inequality}\nonumber\\
&&
\le C 2^{2h}(\sum_{s=h+\ell}^{\infty}2^{-s})\left(\sum_{s=h+\ell}^{\infty}2^{-s}(\int_{2^s\le|x|\le 2^{s+1}} | g(x)|^2dx\right)\nonumber\\
&&
\le
C  2^{h-\ell}\left(\sum_{s=h+\ell}^{\infty}2^{-s}\int_{2^s\le|x|\le 2^{s+1}} | g(x)|^2dx\right)\,.\nonumber
\end{eqnarray}
We observe that in \rec{estrot4bbis} we use the fact that, since $\ell\ge 3$ then $|x-t|,|x-s|\ge 2^{s-1}$ for
every $x,y\in A_h$ and $2^s\le |\xi|\le 2^{s+1}\,.$\par

We multiply the last term in \rec{estrot4bbis} by $2^{k-h}$, where $k\in \Z$,    and we 
 sum  up from $h=k$ to $+\infty$. We get
 \begin{eqnarray}\label{estrot4bis}
&&
\sum_{h=k}^{+\infty}2^{k-h}2^{h-\ell}\left(\sum_{s=h+\ell}^{\infty}2^{-s}\int_{2^s\le|x|\le 2^{s+1}}| M\Delta^{1/4} u|^2dx\right)\nonumber\\
&&
=2^{-\ell}\sum_{s=k+\ell}^{+\infty}2^{k-s}(s-\ell-k)\left(\int_{2^s\le|x|\le 2^{s+1}}| M\Delta^{1/4} u|^2dx\right)\\
&&
\le C2^{-\ell}\sum_{s=k+\ell}^{+\infty}2^{\frac{k-s}{2}}\left (\int_{2^s\le|x|\le 2^{s+1}}| M\Delta^{1/4} u|^2dx\right)\,.\nonumber
\end{eqnarray}
 
{\bf 2.}  Estimate of  $\int_{A_h}\int_{A_h}\frac{1}{|t-s|^2}\left(\int_{|x|<2^{h-\ell}}  g(x)(\frac{1}{|x-t|^{1/2}}-\frac{1}{|x-s|^{1/2}})dx\right)^2  dt ds$\,.\par
 For $h\ge k$  we have
\begin{eqnarray}\label{estrot4tris}
&&
\int_{A_h}\int_{A_h}\frac{1}{|t-s|^2}\left(\int_{|x|<2^{h-\ell}}  g(x)(\frac{1}{|x-s|^{1/2}}-\frac{1}{|x-t|^{1/2}})dx\right)^2  dt ds\nonumber\\
&&\mbox{ny Mean Value Theorem}\nonumber\\&&
\le C \int_{A_h}\int_{A_h}\left(\int_{|x|<2^{h-\ell}}  g(x)\max(\frac{1}{|x-t|^{3/2}},\frac{1}{|x-s|^{3/2}})dx\right)^2  dt ds  \\
&&\le
C \int_{A_h}\int_{A_h} 2^{-3h}2^{h-\ell}\left(\int_{|x|<2^{h-\ell}} | g(x)|^2 dx\right)dt ds\nonumber\\&&
=
C2^{-\ell}  \int_{|x|<2^{h-\ell}}| M\Delta^{1/4} u|^2 dx\nonumber\\
&&=
C2^{-\ell} \left( \int_{|x|<2^{k-\ell}}| M\Delta^{1/4} u|^2 dx+\sum_{s=k-\ell}^{h-\ell} \int_{2^{s}\le|\xi|<2^{s+1}}| M\Delta^{1/4} u|^2 dx\right)\,.\nonumber
\end{eqnarray}
In \rec{estrot4tris} we use the fact that since $\ell\ge 3$, $t,s\in A_h$ and $|x|<2^{h-\ell}$ we have $|x-s|,|x-t|\ge 2^{h-2}\,.$\par
We multiply   \rec{estrot4tris} by $2^{k-h}$, and we 
 sum  up from $h=k$ to $+\infty$. We get\par
  \begin{eqnarray}\label{estrot4ttris}
  &&\int_{A_h}\int_{A_h}\frac{1}{|x-y|^2}\left(\int_{|x|<2^{h-\ell}}  g(x)(\frac{1}{|x-s|^{1/2}}-\frac{1}{|x-t|^{1/2}})dx\right)^2  dt ds\nonumber\\
&&\le C2^{-\ell}
 \int_{|x|<2^{k-\ell}}| M\Delta^{1/4} u|^2dx+ C2^{-2\ell}\sum_{h=k-\ell}^{+\infty}2^{k-h} \int_{2^{h}\le|x|\le 2^{h+1}} | M\Delta^{1/4} u|^2 dx\,.\
\end{eqnarray}
We choose  $\ell$   so  that  $C2^{-\ell}<\gamma$\, and let $ n_0\ge \ell$. Then  for all $n\ge  n_0$ we obtain
\begin{eqnarray*}
  &&\sum_{h=k}^{+\infty}2^{k-h}
\left[C2^{-\ell}\int_{|x|<2^{k-\ell}}| M\Delta^{1/4} u|^2dx+ C2^{-2\ell}\sum_{s=k-\ell}^{h-\ell} \int_{2^{s}\le|x|\le 2^{s+1}}| M\Delta^{1/4} u|^2 dx\right]\\
&&\le \gamma
 \int_{|x|<2^{k-n}}| M\Delta^{1/4} u|^2dx+ \sum_{h=k-n}^{+\infty} 2^{k-h}\int_{2^{h}\le|x|\le 2^{h+1}} | M\Delta^{1/4} u|^2 dx\,.\
\end{eqnarray*}
By combining \rec{estrDbis}, \rec{estrot4bis}, \rec{estrot4ttris}, for $n\ge n_0$   we finally get 
 \begin{eqnarray*}
&&\sum_{h=k}^{+\infty} 2^{k-h}||\Delta^{-1/4}( M\Delta^{1/4} u)||^2_{\dot{H}^{1/2}(A_h)}\\
&&\le\gamma \int_{|x|<2^{k-n}}| M\Delta^{1/4} u|^2dx+  \sum_{h=k-n}^{+\infty} \int_{2^{h-1}\le|x|\le 2^{h+1}}2^{k-h}| M\Delta^{1/4} u|^2 dx\,.\
\end{eqnarray*}
 and  we conclude the proof.~~\hfill$\Box$
\par\medskip
   Next we show a sort of Poincar\'e Inequality for functions in $\dot{H}^{1/2}(\R)$ having compact support.
   We remark that in general the extension by zero of a function in  $H_0^{1/2}(\Omega)=\overline{C_0^\infty(\Omega)}^{H^{1/2}}$,
    $\Omega$ open subset of $\R$ is not in $H^{1/2}(\R)\,.$ This is the reason why Lions and Magenes \cite{LM} introduced the set $H^{1/2}_{00}(\Omega)$ for which Poincar\'e Inequality holds.
   \begin{Theorem}\label{poincare}
   Let $v\in \dot{H}^{1/2}(\R)$ be such that $\mbox{supp$(v)$}\subset (-1,1)$.\par\    Then $v\in L^2([-1,1])$ and
   $$
   \int_{[-1,1]} |v(x)|^2 dx\le C ||v||^2_{\dot{H}^{1/2}}((-2,2))\,.
   $$
   \end{Theorem}
   {\noindent \bf  Proof.}
   We have
   \begin{eqnarray*}
   &&
     \int_{[-1,1]} |v(x)|^2 dx\le 9\int_{1\le |y|\le 2|}\int_{|x|\le 1} \frac{|v(x)|^2}{|x-y|^2} dx dy\\
     &&\le C \int_{1\le |y|\le 2}\int_{|x|\le 1} \frac{|v(x)|^2}{|x-y|^2} dx dy\\
     && \le C\int_{1\le |y|\le 2|}\int_{|x|\le 1} \frac{|v(x)-v(y)|^2}{|x-y|^2} dx dy\\
     && \le C \int_{ |y|\le 2|}\int_{|x|\le 2} \frac{|v(x)-v(y)|^2}{|x-y|^2} dx dy=C ||v||^2_{\dot{H}^{1/2}}([-2,2])\,.
     \end{eqnarray*}
     We can conclude.~\hfill $\Box$\par
     
     From Lemma \ref{poincare}  it follows that 
     $$
      ||v||_{L^2((-r,r))}\le C r^{1/2}||v||_{\dot{H}^{1/2}(\R)} \,.
      $$
   
       We conclude this Section with the following technical result.
   \begin{Proposition}\label{seq}
   Let $(a_k)_k$ be a sequence  of positive real  numbers satisfying $\sum_{k=-\infty}^{+\infty} a^2_k<+\infty$
   and for every $n\le 0$ 
   \begin{equation}\label{sequence}
    \sum_{-\infty}^{n}a^2_k \le C\left(\sum_{k=n+1}^{+\infty}2^{\frac{n+1-k}{2}} a^2_k\right)\,.
   \end{equation}
   Then there are  $0<\beta<1$ , $C>0$ and $\bar n<0$ such that for $n\le \bar n$ we have  
   $$
    \sum_{-\infty}^{n}a^2_k \le C (2^{n})^{\beta}\,.
    $$
    \end{Proposition}
    {\bf \noindent Proof.}  
    For $n<0$, we set $A_{n}= \sum_{-\infty}^{n}a_k^2$. We have $a^2_k = A_k-A_{k-1}$ and thus
     \begin{eqnarray*}
     &&
     A_{n}\le C\sum_{n+1}^{+\infty}2^{\frac{n+1-k}{2}} (A_k-A_{k-1})\le C(1-1/\sqrt{2})\sum_{n+1}^{+\infty}2^{\frac{n+1-k}{2}} A_k-CA_{n}\,.
     \end{eqnarray*}
     Therefore
     \begin{equation}\label{tau}
     A_{n}\le \tau \sum_{n+1}
     ^{+\infty}2^{\frac{n+1-k}{2}} A_k\,,
     \end{equation}
     $\tau=\frac{C}{(C+1)}(1-1/\sqrt{2}) <1-1/\sqrt{2}\,.$
     \par
     The relation \rec{tau} implies the following estimate
     \begin{eqnarray*}
     &&
    A_{n}\le \tau A_{n+1}+\tau \sum_{n+2}
     ^{+\infty}2^{\frac{n+1-k}{2}} A_k\\
     &&\mbox{by induction}\\
     &&\le 
      \tau^2\left(\sum_{n+2}
     ^{+\infty}2^{\frac{n+2-k}{2}} A_k\right)+\frac{\tau}{\sqrt{2}} \left(\sum_{n+2}
     ^{+\infty}2^{\frac{n+2-k}{2}} A_k\right)\\
     &&
     = \tau(\tau+1/\sqrt{2})\left(\sum_{n+2}
     ^{+\infty}2^{\frac{n+2-k}{2}} A_k\right)\\
     &&
     =\tau(\tau+1/\sqrt{2})\left[A_{n+2}+1/\sqrt{2}\sum_{n+3}
     ^{+\infty}2^{\frac{n+3-k}{2}} A_k\right]
     \\ &&\mbox{again by induction}\\
     &&
     \le \tau(\tau+1/\sqrt{2})^2\sum_{n+3}
     ^{+\infty}2^{\frac{n+3-k}{2}} A_k\\
     &&\le \ldots\\
     &&\le 
     \tau(\tau+1/\sqrt{2})^{-n}\sum_{k=0}^{+\infty}2^{-k} A_k
     \\
     &&
     \le \tau(\tau+1/\sqrt{2})^{-n}
     \sum_{h=-\infty}^{+\infty} a^2_h\,.
     \end{eqnarray*}
     
           Therefore  for some $\beta\in (0,1)$ and for all $n<0$\ we have
     $$ A_{n}\le C (2^{n})^{\beta}\,.~~~~~~\hfill\Box$$
     

\section{$L$-Energy Decrease Controls.}\label{regulsection}

In this Section we provide some {\em localization estimates} of solutions to the following equations
  \begin{equation}\label{eq1intrbis}
\Delta^{1/4} ( M\Delta^{1/4}u)=T(Q,u)\,;
\end{equation}
and
\begin{equation} \label{eq2intrbis}
\Delta^{1/4} (M \Delta^{1/4}u)=S(Q,u)-{\cal{R}}(\Delta^{1/4} u \cdot {\cal{R}}\Delta^{1/4} u)\,,
\end{equation}
where $Q\in \dot H^{1/2}(\R,{\cal{M}}_{\ell\times m}(\R)) $ , $M\in \dot H^{1/2}(\R,{\cal{M}}_{p\times m}(\R)),$ $\ell,p\ge 1$\,. 
\par\medskip
We will consider a dyadic decomposition of the unity   $\varphi_j$    such that
  $$\mbox{ supp$(\varphi_j) \subset B_{2^{j+1}}\setminus B_{2^{j-1}}$},~~\sum_{-\infty}^{+\infty}\varphi_j=1\,.$$ \par
   For every $k,h\in\Z$, we set
    $$\chi_k:=\sum_{-\infty}^{k-1}\varphi_j\,,,~\bar{u}_k=|B_{2^k}|^{-1}\int_{B_{2^k}} u(x)\, dx,$$
    $$A_h=B_{2^{h+1}}\setminus B_{2^{h-1}}
  ~\mbox{ and}~ \bar{u}^{h}=|A_h|^{-1}\int_{A_h} u(x)\, dx\,,$$
$$
A^\prime_h=B_{2^{h}}\setminus B_{2^{h-1}}
  ~\mbox{ and}~ \bar{u}^{\prime,h}=|A^\prime_h|^{-1}\int_{A^\prime_h} u(x)dx\,.$$
   \medskip
We prove the following results.
    
\begin{Lemma}\label{regularity1}
Let  $Q\in \dot H^{1/2}(\R,{\cal{M}}_{\ell\times m}(\R))$, $M\in \dot H^{1/2}(\R,{\cal{M}}_{p\times m}(\R)) ,$  $\ell,p\ge 0$ and let $u\in \dot H^{1/2}(\R,\R^m)$ be a solution
of \rec{eq1intrbis}. Then for  $k<0$ with $|k|$  large enough  we  have
\begin{eqnarray}\label{estrot5}
 ||M\Delta^{1/4} u||^2_{L^2 (B_{2^k})}-\frac{1}{4}||\Delta^{1/4} u||^2_{L^2 (B_{2^k})}&\le &C\left[\sum_{h=k}^\infty (2^\frac{k-h}{2})||M\Delta^{1/4} u||^2_{L^2(A_h)}\right. \\
&& \left.+\sum_{h=k}^\infty (2^\frac{k-h}{2})||\Delta^{1/4} u||^2_{L^2(A_h)}\right]\,.\nonumber\end{eqnarray}
\end{Lemma}
 \begin{Lemma}\label{regularity2}
Let  $Q\in \dot H^{1/2}(\R,{\cal{M}}_{\ell\times m})$, $M\in \dot H^{1/2}(\R,{\cal{M}}_{p\times m}(\R)),$  $\ell,p\ge 1$ and let $u\in \dot H^{1/2}(\R,\R^m)$ be a solution
of \rec{eq2intrbis}. Then for  $k<0$ with $|k|$  large enough  the estimate \rec{estrot5} holds.
 \end{Lemma}
In the next Section we will see  that weak $1/2$-harmonic maps $u$ satisfy both the  equations \rec{eq1intr} and \rec{eq2intr} which are \rec{eq1intrbis} and \rec{eq2intrbis}
with $(M,Q)$ given respectively by
 $(u\wedge,u\wedge)$ and $(u\cdot,u\cdot)$\,.   \medskip

   We premise some estimates.
   \begin{Lemma}\label{estphi}
   Let $u\in \dot H^{1/2}(\R)$. Then for all $k\in\Z$
  the following estimate holds  \begin{equation}\label{phih1}
 \sum_{h=k}^{+\infty} 2^{k-h} ||\varphi_h(u-\bar u_k)||_{\dot{H}^{1/2}(\R)}\le C\left[\sum_{s\le k } 2^{s-k}||u||_{\dot{H}^{1/2}(A_s)}+ \sum_{s\ge k } 2^{k-s}||u||_{\dot{H}^{1/2}(A_s)}\right]\,.\end{equation}
   \end{Lemma}
{\bf Proof of Lemma \ref{estphi}\,.}
We have first
\begin{eqnarray}\label{phih0}
&&  ||\varphi_h(u-\bar u_k)||_{\dot{H}^{1/2}(\R)}\le ||\varphi_h(u-\bar u^h)||_{\dot{H}^{1/2}(\R)}+
||\varphi_h||_{\dot{H}^{1/2}(\R)}|\bar  u_k-\bar u^h|\,.
\end{eqnarray}
 We estimate the r.h.s of \rec{phih0}\,.
We have
\begin{eqnarray}\label{phih}
&& ||\varphi_h(u-\bar u^h)||_{\dot{H}^{1/2}(\R)}\\
&&=\int_{A_h}\int_{A_h}\frac{|\varphi_h(u-\bar u^h)(x)-\varphi_h(u-\bar u^h)(y)|^2}{|x-y|^2} dx dy\nonumber\\
&&\le
2 \left[\int_{A_h}\int_{A_h}\frac{|u(x)-u(y)|^2}{|x-y|^2} dx dy+||\nabla \varphi_h||^2_{\infty}\int_{A_h}\int_{A_h} |u-\bar u^h|^2 dx dy\right]\nonumber\\\
&& \le C \left [||u||^2_{\dot{H}^{1/2}(A_h)}+2^{-h}\int_{A_h} |u-\bar u^h|^2 dx\right]\nonumber\\
&&\le C ||u||^2_{\dot{H}^{1/2}(A_h)}\,.\nonumber\
\end{eqnarray}
where we use the fact $||\nabla \varphi_h||_{\infty}\le C2^{-h}$. \par
Now we estimate $|\bar  u_k-\bar u^h|$.   
We  can write
$\bar {u}_k=\sum_{\ell=-\infty}^{k-1}2^{\ell-k}\bar{u}^{\prime,\ell}\,.$ \par
Moreover
\begin{eqnarray}\label{phihbis}
&& |\bar  u_k-\bar u^h|\le |\bar  u^h-\bar u^{\prime,h}|+|\bar  u_k-\bar u^{\prime,h}|\nonumber\\
&&\le C|A_h|^{-1}\int_{A_h}|u-\bar u^h|\, dx+  \sum_{\ell=-\infty}^{k-1} 2^{\ell-k}\sum_{s=\ell}^{h-1}|\bar u^{\prime,s+1}-\bar u^{\prime,s}|\nonumber\\
&&\le C|A_h|^{-1}\int_{A_h}|u-\bar u^h|\, dx+\sum_{\ell=-\infty}^{k-1} 2^{\ell-k}\sum_{s=\ell}^{h-1} |A_{s+1}|^{-1}\int_{A_{s+1}}|u-\bar{u}^{s+1}|\,dx\\
&& \le C\left[ ||u||_{\dot{H}^{1/2}(A_h)}+\sum_{\ell=-\infty}^{k-1} 2^{\ell-k}\sum_{s=\ell}^{h-1}||u||_{\dot{H}^{1/2}(A_{s+1})}\right]\,.\nonumber
\end{eqnarray}
 Thus combining \rec{phih} and \rec{phihbis} we get 
 \begin{eqnarray}\label{suppne}
 &&  ||\varphi_h(u-\bar u^h)||_{\dot{H}^{1/2}(\R)}\le  [||\varphi_h(u-\bar u^h)||_{\dot{H}^{1/2}(\R)}+
||\varphi_h||_{\dot{H}^{1/2}(\R)}|\bar  u_k-\bar u^h|]\nonumber \\
&& \le C\left[ ||u||_{\dot{H}^{1/2}(A_h)}+\sum_{\ell=-\infty}^{k-1} 2^{\ell-k}\sum_{s=\ell}^{h-1}||u||_{\dot{H}^{1/2}(A_{s+1})}\right]\,.
\end{eqnarray}
Multiplying both sides of \rec{suppne} by $2^{k-h}$ and  summing uo from $h=k$ to $+\infty$ we get
\begin{eqnarray}\label{sumphih}
&&  \sum_{h=k}^{+\infty}2^{k-h}\sum_{\ell=-\infty}^{k-1} 2^{\ell-k}\sum_{s=\ell+1}^{h}||u||_{\dot{H}^{1/2}(A_s)}\\
&& \le C \sum_{s\le k }||u||_{\dot{H}^{1/2}(A_s)} \sum_{h\ge k}  \sum_{\ell\le s} 2^{\ell-h}+ \sum_{s\ge k }||u||_{\dot{H}^{1/2}(A_s)} \sum_{h\ge s }  \sum_{\ell\le k} 2^{\ell-h}\nonumber
\\
&&\le C \sum_{s\le k } 2^{s-k}||u||_{\dot{H}^{1/2}(A_s)}+ \sum_{s\ge k } 2^{k-s}||u||_{\dot{H}^{1/2}(A_s)}\,.\nonumber
\end{eqnarray}
This ends the proof of Lemma  \ref{estphi}.\,.~~\hfill$\Box$

\par\bigskip
 Now we  recall the value of the  Fourier transform of some functions that will be used in the sequel.\par
We have 
${\cal{F}}[|x|^{-1/2}](\xi)=|\xi|^{-1/2}$.
The Fourier transforms of $|x|$, $x|x|^{-1/2}$, $|x|^{1/2}$ are the tempered distributions  defined, for every $\varphi\in {\cal{S}}(\R)$, respectively  by
 \begin{eqnarray}\label{conv}
\langle {\cal{F}}[|x| ],\varphi\rangle&=&
\langle{\cal{F}} [\frac{x}{|x|}]\star{\cal{F}} [x],\varphi\rangle=\langle p.v. (\frac{1}{x})\star({\delta})^\prime_0(x),\varphi\rangle\\ &=&
p.v. \int_{\R}\frac{\varphi(x)-\varphi(0)}{x^2}dx\,;\nonumber
\end{eqnarray}
\begin{eqnarray}\label{conv2}
\langle {\cal{F}}[x|x|^{-1/2} ],\varphi\rangle &=&
\langle{\cal{F}} [{x}]\star{\cal{F}} [|x|^{-1/2}],\varphi\rangle=\langle p.v.(\frac{1}{x})\star({\delta})^\prime_0(x),\varphi\rangle\\ &=&
p.v.\int_{\R}[\varphi(x)-\varphi(0)]{\frac{x}{|x|}}\frac{1}{|x|^{3/2}}dx\,
\end{eqnarray}
and
 $$\langle {\cal{F}}[|x|^{1/2}],\varphi\rangle=p.v.\int_{\R}\frac{\varphi(x)-\varphi(0)}{|x|^{3/2}} dx\,.$$
  
 \par
\bigskip

Next we set
 $$ F(Q,a)=\Delta^{1/4}(Qa)-Q\Delta^{1/4} a+\Delta^{1/4} Q a\,,$$
 and
 $$
 G(Q,a)={\cal{R}}\Delta^{1/4}(Qa)-Q\Delta^{1/4}{\cal{R}}a+ \Delta^{1/4} Q {\cal{R}}a\,.$$\par\medskip
 We observe that  $T(Q,u)=F(Q,\Delta^{1/4} u)$ and 
 $S(Q,u)={\cal{R}}G(Q,\Delta^{1/4} u)\,.$\par 
\begin{Lemma}\label{reg1}
Let $Q\in \dot{H}^{1/2}(\R,{\cal{M}}_{\ell\times m}(\R))$, $M\in \dot{H}^{1/2}(\R,{\cal{M}}_{p\times m}(\R))$ and let $u\in \dot H^{1/2}(\R)$ be a solution
of \rec{eq1intrbis}.
 Then there exist $C>0$, $\bar n>0$ (independent of $u$ and $M$) such that  
 for all $\eta\in (0,1/4)$ for all $k<k_0$ ($k_0$ depending on $\eta$)   and $n\ge \bar n\,,$  we have
\begin{eqnarray}\label{esttimes}
&&||\chi_{k-4}(w-\bar w_{k-4})||_{\dot{H}^{1/2}(\R)} \le \eta ||\chi_{k-4}\Delta^{1/4}u ||_{L^2}\\
&&+C
\left(\sum_{h=k}^\infty 2^\frac{k-h}{2}||\Delta^{1/4} u||_{L^2(A_h)}+\sum_{h=k-n}^{+\infty}2^{k-h}||w||_{\dot{H}^{1/2}(A_h)}\right)\nonumber
\end{eqnarray}
where $w=\Delta^{-1/4}( M \Delta^{1/4} u)$\, and we recall that $\chi_{k-4}\equiv 1$ on $B_{2^{k-5}}$ and $\chi_{k-4}\equiv 0$
on $B_{2^{k-4}}^c$\,.
\end{Lemma}\par\medskip
 \begin{Lemma}\label{reg2}
Let $Q\in \dot{H}^{1/2}(\R,{\cal{M}}_{\ell\times m}(\R))$, $M\in \dot{H}^{1/2}(\R,{\cal{M}}_{p\times m}(\R))$ and let $u\in \dot{H}^{1/2}(\R)$ be a solution
of \rec{eq2intrbis}. Then there exist $C>0$, $\bar n>0$ (independent of $u$ and $M$)  such that  
 for all $\eta\in (0,1/4)$, for all $k<k_0$ ($k_0$ depending on $\eta$)   and $n\ge \bar n\,,$  we have
\begin{eqnarray}\label{estdot}
&&||\chi_{k-4}(w-\bar w_{k-4})||_{\dot{H}^{1/2}(\R)} \le \eta||\chi_{k-4}\Delta^{1/4}u ||_{L^2(\R)}\\
&&+ C
\left(\sum_{h=k}^\infty (2^\frac{k-h}{2} ||\Delta^{1/4} u||_{L^2(A_h)}+\sum_{h=k-n}^{k-3}2^{h-k}||w||_{\dot{H}^{1/2}(A_h)}\right)\nonumber
\end{eqnarray}
where $w=\Delta^{-1/4}( M \Delta^{1/4} u)$\,.
\end{Lemma}
{\bf Proof of Lemma  \ref{reg1}\,.} \par
We fix $\eta\in (0,1/4)$. \par
 Let $k<0$ be large enough so that
 $||\chi_k(Q-\bar Q_k)||_{\dot{H}^{1/2}(\R)}\le \varepsilon$, where $\varepsilon\in (0,1)$ will be determined later.
 \par
   We write
 $$F(Q,\Delta^{1/4}u)=F(Q_1,\Delta^{1/4}u)+F(Q_2,\Delta^{1/4}u)\,,$$
 where
 $Q_1=\chi_k(Q-\bar Q_k)$ and $Q_2=(1-\chi_k)(Q-\bar Q_k)$\,. We observe that, by construction, we have 
 $supp( Q_2)\subseteq B^c_{2^{k-1}}$, $||Q_1||_{\dot{H}^{1/2}(\R)}\le \varepsilon\,$ and  $||Q_2||_{\dot{H}^{1/2}(\R)}\le ||Q||_{\dot{H}^{1/2}(\R)}\,.$ \par
  We rewrite the equation \rec{eq1intrbis}
 as follows:
 \begin{eqnarray}\label{eq3}
 && \Delta^{1/2}(\chi_{k-4} (w-\bar w_{k-4}))=-\sum_{h=k-4}^{+\infty} \Delta^{1/2}(\varphi_h (w-\bar w_{k-4}))\\
 && ~~~~~+F(Q_1,\Delta^{1/4} u)+F(Q_2,\Delta^{1/4} u)\,.\nonumber
 \end{eqnarray}
 We multiply the equation \rec{eq3} by $\chi_{k-4} (w-\bar w_{k-4})$ and integrate over $\R$\,.
 We get
 \begin{eqnarray}\label{eqint1}
 &&
 \int_{\R}|\Delta^{1/4}(\chi_{k-4} (w-\bar w_{k-4}))|^2dx =
 -\sum_{h=k-4}^{+\infty}\int_{\R} \Delta^{1/2}(\varphi_h (w-\bar w_{k-4}))(\chi_{k-4} (w-\bar w_{k-4}))dx\nonumber\\
 &&
 +\int_{\R}F(Q_1,\Delta^{1/4} u)(\chi_{k-4} (w-\bar w_{k-4}))dx+  \int_{\R} F(Q_2,\Delta^{1/4} u)(\chi_{k-4} (w-\bar w_{k-4})) dx\,.
 \end{eqnarray}
 We estimate the last three terms     in \rec{eqint1}\,.
 \par
 $\bullet$ { Estimate of $-\sum_{h=k-4}^{+\infty}\int_{\R} \Delta^{1/2}(\varphi_h (w-\bar w_{k-4}))(\chi_{k-4} (w-\bar w_{k-4}))$\,.}\par
 {\bf Case $k-4\le h\le k-3$\,.}
 \begin{eqnarray}\label{estintphi}
 &&\sum_{h=k-4}^{k-3}\int_{\R} \Delta^{1/2}(\varphi_h (w-\bar w_{k-4}))(\chi_{k-4} (w-\bar w_{k-4}))dx \nonumber\\
 &&
 \le ||(\chi_{k-4} (w-\bar w_{k-4}))||_{\dot{H}^{1/2}(\R)} \left[
  \sum_{h=k-4}^{k-3}  || (\varphi_h (w-\bar w_{k-4}))||_{\dot{H}^{1/2}(\R)} \right] 
\nonumber  \\
  &&\mbox{by 
Lemma \ref{estphi}}
\\
&& \le || (\chi_{k-4} (w-\bar w_{k-4}))||_{\dot{H}^{1/2}(\R)} \left(
  \sum_{h=k-4}^{k-3}\left[||w||_{\dot{H}^{1/2}(A_h)}+\sum_{\ell=-\infty}^{k-5} 2^{\ell-(k-4)}\sum_{s=\ell+1}^{h}||w||_{\dot{H}^{1/2}(A_s)}\right]\right)\nonumber \\
&&\le C|| (\chi_{k-4} (w-\bar w_{k-4}))||_{\dot{H}^{1/2}(\R)} \left[\sum_{h=-\infty}^{k-3} 2^{h-k}||w||_{\dot{H}^{1/2}(A_h)}\right] \nonumber
\end{eqnarray}
From Localization Theorem \ref{localization} it follows that
$$
\sum_{h=-
\infty}^{k-6}  ||w||^2_{\dot{H}^{1/2}(A_h)} \le C|| (\chi_{k-4} (w-\bar w_{k-4}))||^2_{\dot{H}^{1/2}(\R)}\,,$$
where $C>0$ is independent of $k$ and $w\,.$
Thus there exists 
  $n_1\ge 6$  such that  $n\ge   n_1$  we have
$$
C  \sum_{h=-
\infty}^{k-n} 2^{h-k}||w||_{\dot{H}^{1/2}(A_h)} \le \frac{1}{8}|| (\chi_{k-4} (w-\bar w_{k-4}))||_{\dot{H}^{1/2}(\R)}\,.
$$
Thus for  $n\ge   n_1$ we have
\begin{eqnarray}\label{estintphibis}
&&\sum_{h=k-4}^{k-3}\int_{\R} \Delta^{1/2}(\varphi_h (w-\bar w_{k-4}))(\chi_{k-4} (w-\bar w_{k-4}))\nonumber\\
 &&\le \frac{1}{8}|| (\chi_{k-4} (w-\bar w_{k-4}))||^2_{\dot{H}^{1/2}(\R)}+C \sum_{h=
k-n}^{k-3} 2^{h-k}||w||_{\dot{H}^{1/2}(A_h)}\,.
\end{eqnarray}
{\bf Case $k-2\le h<+\infty$\,.}\par
 In this case we use the fact that the supports of $\varphi_h$ and of $\chi_{k-4}$ are disjoint and in particular
 $0\notin (\varphi_h (w-\bar w_{k-4}))\star (\chi_{k-4} (w-\bar w_{k-4}))$.\par
\begin{eqnarray}\label{estintphi2}
&&\sum_{h=k-2}^{+\infty}\int_{\R} \Delta^{1/2}(\varphi_h (w-\bar w_{k-4}))(\chi_{k-4} (w-\bar w_{k-4}))dx \nonumber\\
&&
\sum_{h=k-2}^{+\infty}\int_{\R}{\cal{F}}^{-1}(|\xi|)(x)(\varphi_h (w-\bar w_{k-4}))\star (\chi_{k-4} (w-\bar w_{k-4}))dx\nonumber\\
&&
\le 
\sum_{h=k-2}^{+\infty}||{\cal{F}}^{-1}(|\xi|)||_{L^{\infty}(B_{2^{h+2}}\setminus B_{2^{h-2}})}||\varphi_h (w-\bar w_{k-4})||_{L^1}||\chi_{k-4} (w-\bar w_{k-4}|)|_{L^{1}}\nonumber\\
&&
\le C\sum_{h=k-2}^{+\infty}2^{-2h}2^{h/2}||\varphi_h (w-\bar w_{k-4})||_{L^2(\R)}2^{k/2}||\chi_{k-4} (w-\bar w_{k-4})||_{L^2(\R)}\nonumber\\&& 
\mbox{by Theorem \ref{poincare}}\\&&
\le 
 C \sum_{h=k-2}^{+\infty} 2^{k-h}||\varphi_h (w-\bar w_{k-4}||_{\dot{H}^{1/2}(\R)} ||\chi_{k-4} (w-\bar w_{k-4}||_{\dot{H}^{1/2}(\R)}\nonumber
\\
 &&\mbox{by 
Lemma \ref{estphi}}
\nonumber\\
&&
\le C  \sum_{h=k-2}^{+\infty} 2^{k-4-h}\left[||w||_{\dot{H}^{1/2}(A_h)}+\sum_{\ell=-\infty}^{k-5} 2^{\ell-(k-4)}\sum_{s=\ell+1}^{h}||w||_{\dot{H}^{1/2}(A_s)}\right]\nonumber\\
&&~~~~~~~||\chi_{k-4} (w-\bar w_{k-4})||_{\dot{H}^{1/2}(\R)} \nonumber
\\
&&
\le C  \left[ \sum_{h=k-2}^{+\infty} 2^{k-4-h}||w||_{\dot{H}^{1/2}(A_h)}+\sum_{s\le k-4}||w||_{\dot{H}^{1/2}(A_s)}\left(\sum_{h\ge k-4}\sum_{\ell\le s-1} 2^{\ell-h}\right)\right.\nonumber\\
&&
\left.+\sum_{s\ge k-4}||w||_{\dot{H}^{1/2}(A_s)}\left(\sum_{h\ge s-1}\sum_{\ell\le k-4} 2^{\ell-h}\right)\right]||\chi_{k-4} (w-\bar w_{k-4})||_{\dot{H}^{1/2}(\R)} \nonumber\\
&&
\le C \left[ \sum_{h=k-4}^{+\infty} 2^{k-4-h}||w||_{\dot{H}^{1/2}(A_h)}+  \sum_{h=-\infty}^{k-5} 2^{h-(k-4)}||w||_{\dot{H}^{1/2}(A_h)}\right]||\chi_{k-4} (w-\bar w_{k-4})||_{\dot{H}^{1/2}(\R)} \nonumber
\end{eqnarray}
Let $ n_2\ge 6 $ be such that   if $n\ge n_2$ we have 
$$
C \sum_{h=-\infty}^{k- n} 2^{h-(k-4))}||w||_{\dot{H}^{1/2}(A_h)}\le \frac{1}{8} ||\chi_{k-4} (w-\bar w_{k-4})||_{\dot{H}^{1/2}(\R)}\,.$$
 Thus if $n>\bar n=\max({n_1,n_2})$, then from \rec{estintphibis} and \rec{estintphi2} it follows  
 
\begin{eqnarray}\label{estintphi3}
&&
\sum_{h=k-4}^{+\infty}\int_{\R} \Delta^{1/2}(\varphi_h (w-\bar w_{k-4}))(\chi_{k-4} (w-\bar w_{k-4}))\\
&&\le \frac{1}{4} ||\chi_{k-4} (w-\bar w_{k-4})||^2_{\dot{H}^{1/2}(\R)} +C\sum_{h=k-n}^{+\infty} 2^{k-h}||w||_{\dot{H}^{1/2}(A_h)}\nonumber\,.
\end{eqnarray}

$\bullet$ {Estimate of $\int_{\R}F(Q_1,\Delta^{1/4} u)(\chi_{k-4} (w-\bar w_{k-4}))dx$\,.}\par
We write
$$
F(Q_1,\Delta^{1/4} u)=F(Q_1,\chi_{k-4}\Delta^{1/4} u)+\sum_{h=k-4}^{+\infty}F(Q_1,\varphi_h \Delta^{1/4} u)\,.$$
The following estimate holds
\begin{eqnarray}\label{estF1chi}
&&\int_{\R}F(Q_1,\chi_{k-4}\Delta^{1/4} u)(\chi_{k-4} (w-\bar w_{k-4}))dx\\[5mm]
&&\mbox{by Theorem \ref{comm1}}\nonumber\\ 
&&
\le
C||Q_1||_{\dot{H}^{1/2}(\R)}||\chi_{k-4}\Delta^{1/4} u||_{L^2}||\chi_{k-4} (w-\bar w_{k-4})||_{\dot{H}^{1/2}(\R)}\nonumber\\[5mm]
&&  C\varepsilon ||\chi_{k-4}\Delta^{1/4} u||_{L^2}||\chi_{k-4} (w-\bar w_{k-4})||_{\dot{H}^{1/2}(\R)}\,.\nonumber 
\end{eqnarray}
By choosing     $\varepsilon>0$ small enough, we may assume  that  $C\varepsilon<\frac{\eta}{4} <\frac{1}{16}\,.$\par
{\bf Case $k-4\le h\le k+1$.}\par
 We use again Theorem \ref{comm1}.
\begin{eqnarray}\label{estF1phi1}
&&
\sum_{h=k-4}^{k+1}\int_{\R}F(Q_1,\varphi_h \Delta^{1/4} u)(\chi_{k-4}(w-\bar w_{k-4}))dx\\[5mm]
&&\mbox{by Theorem \ref{comm1}}\nonumber\\[5mm] && 
\le
C\sum_{h=k-4}^{k+1} ||Q_1||_{\dot{H}^{1/2}(\R)}||\varphi_h\Delta^{1/4} u||_{L^2(\R)}||\chi_{k-4} (w-\bar w_{k-4})||_{\dot{H}^{1/2}(\R)}\nonumber
\end{eqnarray}\par
{\bf Case $h\ge k+2$.}\par
  We estimate the single terms of  $F(Q_1,\varphi_h \Delta^{1/4} u)(\chi_{k-4}(w-\bar w_{k-4}))\,.$ \par We observe that
if $h\ge k+2$ then the supports of $Q_1$ and $\varphi_h$ and those of $\chi_{k-4}$ and $\varphi_h$ are disjoint. Therefore
$$F(Q_1,\varphi_h\Delta^{1/4} u)\chi_{k-4} (w-\bar w_{k-4})=Q_1\Delta^{1/4}(\varphi_h\Delta^{1/4}u)(\chi_{k-4}(w-\bar w_{k-4}))\,.$$
 
 \begin{eqnarray}\label{estF1phi2}
&&
\sum_{h=k+2}^{+\infty }\int_{\R}F(Q_1,\varphi_h \Delta^{1/4} u)(\chi_{k-4}(w-\bar w_{k-4}))dx\\
&&
=\sum_{h=k+2}^{+\infty }\int_{\R}Q_1\Delta^{1/4}(\varphi_h\Delta^{1/4}u)(\chi_{k-4}(w-\bar w_{k-4}))dx\nonumber\\
&&
=\sum_{h=k+2}^{+\infty }\int_{\R} 
{\cal{F}}^{-1}(|\xi|^{1/2})(x)\left([\varphi_h\Delta^{1/4} u]\star [Q_1(\chi_{k-4}(w-\bar w_{k-4}))]\right)\nonumber\\
&&
\sum_{h=k+2}^{+\infty }||{\cal{F}}^{-1}(|\xi|^{1/2})||_{L^{\infty}(B_{2^{h+2}}\setminus B_{2^{h-2}})}||\varphi_h\Delta^{1/4} u||_{L^1}
||Q_1(\chi_{k-4}(w-\bar w_{k-4}))||_{L^1}\nonumber\\
&&
\le
C\sum_{h=k+2}^{+\infty }2^{-3/2h}||\varphi_h\Delta^{1/4} u||_{L^1}
||Q_1(\chi_{k-4}(w-\bar w_{k-4}))||_{L^1}\nonumber\\
&&\mbox{by Theorem \ref{poincare}}\nonumber\\
&& \le C \sum_{h=k+2}^{+\infty } 2^{k-h}||Q_1||_{\dot{H}^{1/2}(\R)}||\,||\varphi_h\Delta^{1/4} u||_{L^2(\R)}||\chi_{k-4} (w-\bar w_{k-4})||_{\dot{H}^{1/2}(\R)}\nonumber\\
&&
\le
C\sum_{h=k+2}^{+\infty } 2^{k-h}||\varphi_h\Delta^{1/4} u||_{L^2}||\chi_{k-4} (w-\bar w_{k-4})||_{\dot{H}^{1/2}(\R)}\nonumber
\end{eqnarray}

 \par
$\bullet$ { Estimate of $\int_{\R}F(Q_2,\Delta^{1/4} u)\,(\chi_{k-4} (w-\bar w_{k-4}))dx$\,.} \par
As above we write
$$
F(Q_2,\Delta^{1/4} u)=F(Q_2,\chi_{k-4}\Delta^{1/4} u)+\sum_{h=k-4}^{+\infty}F(Q_2,\varphi_h \Delta^{1/4} u)\,.$$

Since the support of $Q_2$ is included in $B^c_{2^{k-1}}$, we have
$$F(Q_2,\chi_{k-4}\Delta^{1/4} u)\,(\chi_{k-4} (w-\bar w_{k-4}))=\Delta^{1/4} Q_2(\chi_{k-4}\Delta^{1/4} u)\,(\chi_{k-4} (w-\bar w_{k-4}))\,.$$
We can write $Q_2=\sum_{h=k-1}^{+\infty}\varphi_h(Q_2-\bar{Q_2}_{k-1})$, ($\bar{Q_2}_{k-1}=0$)\,.
\begin{eqnarray}\label{estF2chi1}
&&
\int_{\R} F(Q_2,\chi_{k-4} \Delta^{1/4} u)(\chi_{k-4}(w-\bar w_{k-4}))dx\\
&&=
\sum_{h=k-1}^{+\infty}\int_{\R} \Delta^{1/4}(\varphi_h(Q_2-\bar{Q_2}_{k-1}))(\chi_{k-4}\Delta^{1/4} u)(\chi_{k-4} (w-\bar w_{k-4}))\nonumber\\
&& 
\le C\sum_{h=k-1}^{+\infty}\int_{\R} {\cal{F}}^{-1}(|\xi|^{1/2})\left([\varphi_h(Q_2-\bar{Q_2}_{k-1})]\star [(\chi_{k-4}\Delta^{1/4} u)(\chi_{k-4} (w-\bar w_{k-4}))]\right)\nonumber
\\[5mm]
&&  
\le C\sum_{h=k-1}^{+\infty}2^{-h/2}2^{k/2}
||\varphi_h(Q_2-\bar{Q_2}_{k-1})||_{\dot{H}^{1/2}(\R)}||\chi_{k-4}\Delta^{1/4} u||_{L^2}||\chi_{k-4} (w-\bar w_{k-4})||_{\dot{H}^{1/2}(\R)}\,.\nonumber
\end{eqnarray}
From Lemma \ref{estphi},  by  choosing possibly   $k$ smaller, it follows that   
$$C\sum_{h=k-1}^{+\infty}2^{\frac{k-h}{2}}||\varphi_h(Q_2-\bar{Q_2}_{k-1})||_{\dot{H}^{1/2}(\R)}\le \frac{\eta}{4}<\frac{1}{16}\,.$$

 $\bullet$ { Estimate of $\sum_{h=k-4}^{+\infty}\int_{\R}F(Q_2,\varphi_h \Delta^{1/4} u)(\chi_{k-4} (w-\bar w_{k-4}))dx$\,.} \par
 The following estimates holds.
 \begin{eqnarray}\label{estF2phi1}
 &&
 \sum_{h=k-4}^{k+1}\int_{\R}F(Q_2,\varphi_h \Delta^{1/4} u)(\chi_{k-4} (w-\bar w_{k-4}))dx\\
 &&
 \le C \sum_{h=k-4}^{k+1}||Q||_{\dot{H}^{1/2}(\R)}||\varphi_h\Delta^{1/4} u||_{L^2}||\chi_{k-4} (w-\bar w_{k-4})||_{\dot{H}^{1/2}(\R)}\,.\nonumber
 \end{eqnarray}
  On another hand  since the support of $Q_2$ is included in $B^c_{2^{k-1}}$, if $h\ge k+1$,  we have
 
 \begin{eqnarray*}
 && F(Q_2,\varphi_h \Delta^{1/4} u)(\chi_{k-4} (w-\bar w_{k-4}))= (\chi_{k-4} (w-\bar w_{k-4}))\\[5mm]
 && \left[\Delta^{1/4}(Q_2\,\varphi_h\Delta^{1/4} u)-Q_2\Delta^{1/4}(\varphi_h\Delta^{1/4} u)+\Delta^{1/4} Q_2\varphi_h\Delta^{1/4} u\right]\\[5mm]&&=
(\chi_{k-4} (w-\bar w_{k-4})) \Delta^{1/4}(Q_2\varphi_h\Delta^{1/4}u) \,.
 \end{eqnarray*}
 Let  $\psi_h\in C^\infty_0(\R)$, $\psi_h\equiv 1$ in $B_{2^{h+1}}\setminus B_{2^{h-1}}$ and
 $supp( \psi) \subset  B_{2^{h+2}}\setminus B_{2^{h-2}}\,.$
  \begin{eqnarray}\label{estF2phi2}
 &&
 \sum_{h=k+1}^{+\infty}\int_{\R}F(Q_2,\varphi_h \Delta^{1/4} u)(\chi_{k-4} (w-\bar w_{k-4}))dx\\
 &&
  \sum_{h=k+1}^{+\infty}\int_{\R} \Delta^{1/4}(Q_2\varphi_h\Delta^{1/4}u) (\chi_{k-4} (w-\bar w_{k-4}))dx\nonumber\\
  &&
  =
 \sum_{h=k+1}^{+\infty}\int_{\R}  {\cal{F}}^{-1}(|\xi|^{1/2})\left([\varphi_h\Delta^{1/4} u(Q_2-\bar{Q_2}_{k-1})]\star [\chi_{k-4} (w-\bar w_{k-4})]\right)dx\nonumber\\
 &&
 \le C \sum_{h=k+1}^{+\infty} 2^{-3/2 h}||\varphi_h\Delta^{1/4} u||_{L^2(\R)}||\psi_h(Q_2-\bar{Q_2}_{k-1})||_{L^2(\R)}||\chi_{k-4} (w-\bar w_{k-4})||_{L^2(\R)}
 \nonumber\\
 && \le
 C\sum_{h=k+1}^{+\infty}2^{k-h}||\varphi^{1/2}_h(Q_2-\bar{Q_2}_{k-1})||_{\dot{H}^{1/2}(\R)}||\varphi_h^{1/2}\Delta^{1/4} u||_{L^2(\R)}||\chi_{k-4} (w-\bar w_{k-4})||_{\dot{H}^{1/2}(\R)}\nonumber\\
 &&
 \mbox{by Cauchy-Schwartz Inequality}\nonumber\\
 &&
\le C\left(\sum_{h=k+1}^{+\infty}2^{k-h}||\psi_h(Q_2-\bar{Q_2}_{k-1})||^2_{\dot{H}^{1/2}(\R)}\right)^{1/2}
\left(\sum_{h=k+1}^{+\infty}2^{k-h}||\varphi_h\Delta^{1/4} u||_{L^2}^2\right)^{1/2}\nonumber\\
&& ~~~~~~~||\chi_{k-4} (w-\bar w_{k-4})||_{\dot{H}^{1/2}(\R)}\,.\nonumber
\end{eqnarray}

From Lemma \ref{estphi} (with $\varphi$ replaced by $\psi$)   and Theorem \ref{localization} we deduce that
 
$$(\sum_{h=k+1}^{+\infty}2^{k-h}||\psi_h(Q_2-\bar{Q_2}_{k-1})||^2_{\dot{H}^{1/2}(\R)})^{1/2}\le C||Q||_{\dot{H}^{1/2}(\R)}\,.$$
Thus
\begin{eqnarray*}
&&
\sum_{h=k+1}^{+\infty}\int_{\R}F(Q_2,\varphi_h \Delta^{1/4} u)(\chi_{k-4} 
(w-\bar w_{k-4}))dx\le C ||\chi_{k-4} (w-\bar w_{k-4})||_{\dot{H}^{1/2}(\R)}\\
&&~~~~~~~\left( \sum_{h=k+1}^{+\infty}2^{k-h}||\varphi_h\Delta^{1/4} u||_{L^2}^2\right)^{1/2}\,.\end{eqnarray*}
By combining \rec{estF1phi1},\rec{estF1phi2}, \rec{estF2phi1} and \rec{estF2phi2}   we obtain (for some constant $C$ depending on $Q$)
\begin{eqnarray}\label{T}
&&\int_{\R}F(Q,\Delta^{1/4} u)(\chi_{k-4} (w-\bar w_{k-4}))dx\le \frac{\eta}{2}| |\chi_{k-4}\Delta^{1/4} u||_{L^2}||\chi_{k-4} (w-\bar w_{k-4})||_{\dot{H}^{1/2}(\R)}\nonumber\\
&& ~~~~~+C\sum_{h=k-4}^{+\infty}2^{\frac{k-h}{2}}|| \Delta^{1/4} u||_{L^2(A_h)}||\chi_{k-4} (w-\bar w_{k-4})||_{\dot{H}^{1/2}(\R)}\,.
\end{eqnarray}
Finally  for all  $n\ge \bar n$ we have 
\begin{eqnarray}\label{final}
&&
||\chi_{k-4} (w-\bar w_{k-4})||_{\dot{H}^{1/2}(\R)}\le \eta||\chi_{k-4}\Delta^{1/4} u||_{\dot{H}^{1/2}(\R)}+ C \sum_{h=k-n}^{+\infty} 2^{k-h}||w||_{\dot{H}^{1/2}(A_s)}\\
&&~~~+
 C\sum_{h=k-4}^{+\infty}2^{\frac{k-h}{2}}|| \Delta^{1/4} u||_{L^2(A_h)}\,.\nonumber
 \end{eqnarray}
   and we
 can conclude.~\hfill$\Box$\par
 \bigskip
 We prove Lemma  \ref{reg2}.  \par
 \medskip
 {\bf Proof of Lemma  \ref{reg2}\,.}
 The proof is similar to that of Lemma  \ref{reg1} thus we just sketch it.\par
 We observe that equation \rec{eq2intrbis} is equivalent to
 \begin{equation}\label{eq4}
 {\cal{R}}\Delta^{1/4}(  M\Delta^{1/4}u)=G(Q,\Delta^{1/4}u)-\Delta^{1/4}u \cdot ( {\cal{R}}\Delta^{1/4}u)\,.
 \end{equation}
 We fix $\eta\in (0,1/4)\,.$\par
 Let $k<0$ be such that
 $||\chi_k(Q-\bar Q_k)||_{\dot{H}^{1/2}(\R)}\le \varepsilon$  and  $||\chi_k\Delta^{1/4}u||_{\dot{H}^{1/2}(\R)}\le \varepsilon\,,$
 with $\varepsilon>0$ to be determined later.\par 
 \par
   We write
 $$G(Q,\Delta^{1/4}u)=G(Q_1,\Delta^{1/4}u)+G(Q_2,\Delta^{1/4}u)\,,$$
 where
 $Q_1=\chi_k(Q-\bar Q_k)$ and $Q_2=(1-\chi_k)(Q-\bar Q_k)$\,. We observe that
 $supp( Q_2)\subseteq B^c_{2^{k-1}}$ and $||Q_1||_{\dot{H}^{1/2}(\R)}\le \varepsilon\,.$\par
 We also set $u_1=\chi_k\Delta^{1/4} u$ and $u_2=(1-\chi_k)\Delta^{1/4} u\,$ and  $w=\Delta^{-1/4}(M\Delta^{1/4} u)\,.$\par
 We rewrite the equation \rec{eq4}
 as follows:
 \begin{eqnarray}\label{eq5}
 &&  {\cal{R}}\Delta^{1/2}(\chi_{k-4} (w-\bar w_{k-4}))=-\sum_{h=k-4}^{+\infty}  {\cal{R}}\Delta^{1/2}(\varphi_h (w-\bar w_{k-4}))\\[5mm]
 && +G(Q_1,\Delta^{1/4} u)+G(Q_2,\Delta^{1/4} u)+ u_1( {\cal{R}}\Delta^{1/4}u)+ u_2( {\cal{R}}\Delta^{1/4}u)\,.
 \nonumber
 \end{eqnarray}
 
 We multiply the equation \rec{eq5} by $\chi_{k-4} (w-\bar w_{k-4})$ and integrate over $\R$\,. We get
 \begin{eqnarray}\label{eqintP1}
 &&
 \int_{\R}|\Delta^{1/4}(\chi_{k-4} (w-\bar w_{k-4}))|^2dx \nonumber\\
 &&=
 -\sum_{h=k-4}^{+\infty}\int_{\R} {\cal{R}}\Delta^{1/2}(\varphi_h (w-\bar w_{k-4}))(\chi_{k-4} (w-\bar w_{k-4}))dx\nonumber \\
 &&
 +\int_{\R}G(Q_1,\Delta^{1/4} u)(\chi_{k-4} (w-\bar w_{k-4}))dx\nonumber\\
 &&+  \int_{\R} G(Q_2,\Delta^{1/4} u)(\chi_{k-4} (w-\bar w_{k-4})) dx\\
 && +\int_{\R}u_1( {\cal{R}}\Delta^{1/4}u)(\chi_{k-4} (w-\bar w_{k-4}))dx+\int_{\R}u_2( {\cal{R}}\Delta^{1/4}u)(\chi_{k-4} (w-\bar w_{k-4}))dx\,.\nonumber
 \end{eqnarray}
 We observe that $\int_{\R}u_2( {\cal{R}}\Delta^{1/4}u)(\chi_{k-4} (w-\bar w_{k-4}))dx=0$, having $u_2$ and $\chi_{k-4}$
 supports disjoint\,.\par\medskip
 
 $\bullet$ { Estimate of $-\sum_{h=k-4}^{+\infty}\int_{\R}  {\cal{R}}\Delta^{1/2}(\varphi_h (w-\bar w_{k-4}))(\chi_{k-4} (w-\bar w_{k-4}))dx$\,.}\par\medskip
 {\bf Case $k-4\le h\le k-3$\,.}
 \begin{eqnarray}\label{estintphiP}
 &&\sum_{h=k-4}^{k-3}\int_{\R}   {\cal{R}}\Delta^{1/2}(\varphi_h (w-\bar w_{k-4}))(\chi_{k-4} (w-\bar w_{k-4}))\nonumber\\
 &&
 \le
  \sum_{h=k-4}^{k-3}\int_{\R} ||\Delta^{1/2}(\varphi_h (w-\bar w_{k-4}))||_{\dot H^{-1/2}(\R)}|| (\chi_{k-4} (w-\bar w_{k-4}))||_{\dot{H}^{1/2}(\R)} 
\nonumber  \\
  &&\mbox{by 
Lemma \ref{estphi}}
\\
&&
\le  \sum_{h=k-4}^{k-3}\left[||w||_{\dot{H}^{1/2}(A_h)}+\sum_{\ell=-\infty}^{k-5} 2^{\ell-(k-4)}\sum_{s=\ell+1}^{h}||w||_{\dot{H}^{1/2}(A_s)}\right]\nonumber\\&& ~~
~~~~||(\chi_{k-4} (w-\bar w_{k-4}))||_{\dot{H}^{1/2}(\R)} \nonumber
\end{eqnarray}
Let $n_1\ge 6$ be such that  
$$
C  \sum_{h=-
\infty}^{k-n_1} 2^{h-k}||w||_{\dot{H}^{1/2}(A_h)} \le \frac{1}{8}|| (\chi_{k-4} (w-\bar w_{k-4}))||_{\dot{H}^{1/2}(\R)}\,.
$$
Thus if $n\ge  n_1$ we have
\begin{eqnarray}\label{estintphithis}
 \rec{estintphiP}\le 
 &&\frac{1}{8}|| (\chi_{k-4} (w-\bar w_{k-4}))||^2_{\dot{H}^{1/2}(\R)}+C \sum_{h=
k-n}^{k-3} 2^{h-k}||w||_{\dot{H}^{1/2}(A_h)}\,.
\end{eqnarray}

{\bf Case $k-2\le h<+\infty$\,.}\par
In this case we use the fact that   
$$supp((\varphi_h (w-\bar w_{k-4}))\star (\chi_{k-4} (w-\bar w_{k-4})))\subseteq B_{2^{h+2}}\setminus B_{2^{h-2}}\,,$$
and in particular $0\notin ((\varphi_h (w-\bar w_{k-4}))\star (\chi_{k-4} (w-\bar w_{k-4})))\,.$
\begin{eqnarray}\label{estintphiP2} 
&&\sum_{h=k-2}^{+\infty}\int_{\R}  {\cal{R}}\Delta^{1/2}(\varphi_h (w-\bar w_{k-4}))(\chi_{k-4} (w-\bar w_{k-4}))dx \\
&&
\sum_{h=k-2}^{+\infty}\int_{\R}{\cal{F}}^{-1}(\xi)(x)\left([\varphi_h (w-\bar w_{k-4}))]\star [\chi_{k-4} (w-\bar w_{k-4})]\right)dx\nonumber\\
&&
=\sum_{h=k-2}^{+\infty}\int_{\R}\delta^\prime_0(x)\left(\varphi_h (w-\bar w_{k-4}))\star (\chi_{k-4} (w-\bar w_{k-4}))\right)dx=0\,. \nonumber
\end{eqnarray}

  $\bullet$ { Estimate of $\int_{\R}u_1( {\cal{R}}\Delta^{1/4}u)(\chi_{k-4} (w-\bar w_{k-4}))dx$\,.}\par
  We have
  \begin{eqnarray*}
  &&\int_{\R}u_1( {\cal{R}}\Delta^{1/4}u)(\chi_{k-4} (w-\bar w_{k-4}))dx\\
  &&=\int_{\R}u_1( {\cal{R}}u_1)(\chi_{k-4} (w-\bar w_{k-4}))dx
  +\sum_{h=k}^{+\infty}\int_{\R}u_1( {\cal{R}}\varphi_h\Delta^{1/4}u)(\chi_{k-4} (w-\bar w_{k-4}))dx\,.
  \end{eqnarray*}
  By applying  Lemma \ref{lemmareg} we get
  \begin{eqnarray*}
  &&\int_{\R}u_1( {\cal{R}}u_1)(\chi_{k-4} (w-\bar w_{k-4}))dx\\[5mm]
  && C||u_1( {\cal{R}}u_1)||_{{\cal{H}}}||(\chi_{k-4} (w-\bar w_{k-4}))||_{\dot{H}^{1/2}(\R)}\nonumber\\[5mm]
  &&\le C||u_1||^2_{L^2}||(\chi_{k-4} (w-\bar w_{k-4}))||_{\dot{H}^{1/2}(\R)}\\[5mm]
  &&
  \le C\varepsilon||\chi_{k}\Delta^{1/4}u||_{L^2}||(\chi_{k-4} (w-\bar w_{k-4}))||_{\dot{H}^{1/2}(\R)}\,.
  \end{eqnarray*}
  By choosing $\varepsilon>0$ smaller we may suppose that $C\varepsilon<\frac{\eta}{4}$\,.\par
  Now we observe that for $h\ge k$ the supports of $\varphi_h$ and $\chi_{k-4}$ are disjoint. Thus
  we have
     \begin{eqnarray*}
 &&\sum_{h={k}}^{+\infty} \int_{\R}u_1( {\cal{R}}\varphi_h\Delta^{1/4}u)(\chi_{k-4} (w-\bar w_{k-4}))dx\\
  &&\sum_{h={k}}^{+\infty} \int_{\R}{\cal{F}}^{-1}(\frac{\xi}{|\xi|})\left([\varphi_h\Delta^{1/4}u]\star [u_1(\chi_{k-4} (w-\bar w_{k-4})]\right)dx\\
  &&
  \le C \sum_{h={k}}^{+\infty}|||x|^{-1}||_{L^\infty(B_{2^{h+2}}\setminus B_{2^{h-2}})}  ||\varphi_h\Delta^{1/4}u\star \Delta^{1/4}u_1(\chi_{k-4} (w-\bar
  w_{k-4}))||_{L^1}\\
  &&
  \le C\sum_{h={k}}^{+\infty}2^{-h} 2^{h/2}2^{k/2} ||\varphi_h\Delta^{1/4}u||_{L^2}||u_1||_{L^2}||(\chi_{k-4} (w-\bar
  w_{k-4}))||_{\dot{H}^{1/2}(\R)}\\
  &&
  \le 
  C\varepsilon\sum_{h={k}}^{+\infty} 2^{\frac{k-2}{2}} ||\varphi_h\Delta^{1/4}u||_{L^2} ||(\chi_{k-4} (w-\bar
  w_{k-4}))||_{\dot{H}^{1/2}(\R)}\,.
 \end{eqnarray*}

 The estimate of the terms $$\int_{\R}G(Q_1,\Delta^{1/4} u)(\chi_{k-4} (w-\bar w_{k-4}))dx~\mbox{and}~\int_{\R} G(Q_2,\Delta^{1/4} u)(\chi_{k-4} (w-\bar w_{k-4})) dx$$ are
 analogous of those of $$\int_{\R}F(Q_1,\Delta^{1/4} u)(\chi_{k-4} (w-\bar w_{k-4}))dx~\mbox{and}~\int_{\R} F(Q_2,\Delta^{1/4} u)(\chi_{k-4} (w-\bar w_{k-4})) dx$$
  and so we omit
 them.~~~\hfill$\Box$\par\medskip
 Now we can prove Theorem \ref{regularity1} and \ref{regularity2}\,.\par
\medskip
{\bf Proof of Theorem \ref{regularity1}\,.}
 
From Lemma  \ref{reg1}, it follows that there exist $C>0$ and $\bar n>0$ such that for all  $n>\bar n$, $0<\eta<1/4$, $k<k_0$ ($k_0$ depending on $\eta$),  every solution to \rec{eq1intrbis} satisfies  for some constant $C>0$  \rec{final}
and thus
\begin{eqnarray}\label{finalbis}
&&
||\chi_{k-4} (w-\bar w_{k-4})||^2_{\dot{H}^{1/2}(\R)}\le {\eta^2}||\chi_{k-4}\Delta^{1/4} u||^2_{L^2}+ C \sum_{h=k-n}^{+\infty} 2^{k-h}||w||^2_{\dot{H}^{1/2}(A_h)}\\
&&~~~~+
 C\sum_{h=k-4}^{+\infty}2^{\frac{k-h}{2}}|| \Delta^{1/4} u||^2_{L^2}\,.\nonumber
 \end{eqnarray}
  We fix   $n\ge \bar n$ \par

From Lemma \ref{rotation} it follows that there exist $C_1,C_2>0$ and   $m_1>0$ (independent on $n$, $k$) such that if    $m\ge   m_1$ we have
 \begin{eqnarray}\label{final3}
 && ||\chi_{k-4} (w-\bar w_{k-4})||^2_{\dot{H}^{1/2}(\R)}\ge   \\
   &&~~~~~
   \ge C_1\int_{B_{2^{k-n-m}}}| M\Delta^{1/4} u|^2 dx
    - C_2\sum_{h=k-n-m}^{+\infty}2^{k-h}\int_{B_{2^{h}}\setminus B_{2^{h-1}}}| M\Delta^{1/4} u|^2 dx\,.\nonumber
   \end{eqnarray}
Finally  from Lemma \ref{rotation2} it follows  that there is $C>0$ such that for all $\gamma\in (0,1)$   there exists $m_2>0$ such that if    $m\ge   m_2$       we have 
   \begin{eqnarray}\label{anelli2}
   &&
   \sum_{h=k-n}^{+\infty} 2^{k-h}||w||^2_{\dot{H}^{1/2}(A_h)}=
\sum_{h=k-n}^{+\infty} 2^{k-h}||\Delta^{-1/4}( M\Delta^{1/4} u)||^2_{\dot{H}^{1/2}(A_h)}\\&&
\le \gamma \int_{|\xi|\le 2^{k-n-m}}| M\Delta^{1/4} u|^2 dx
+ \sum_{h=k-n-m}^{+\infty}2^{\frac{k-h}{2}}\int_{2^{h}\le |\xi|\le 2^{h+1}}| M\Delta^{1/4} u|^2 dx\,.\nonumber
\end{eqnarray}
 By   combining \rec{finalbis}, \rec{final3} and \rec{anelli2} we get 
 \begin{eqnarray}\label{estrot7}
C_1|| M\Delta^{1/4} u||^2_{L^2(B_{2^{k-n-m}})}&\le& C\sum_{h=k-n-m}^\infty (2^\frac{k-h}{2})||M\Delta^{1/4} u||^2_{L^2(A_h)}+C_2 \sum_{h=k-n-m}^{+\infty}2^{\frac{k-h}{2}}|| \Delta^{1/4} u||_{L^2(A_h)}\nonumber\\[5mm]
&&+{\eta^2}||\chi_{k-4}\Delta^{1/4} u||^2_{L^2(\R)}+C\gamma|| M\Delta^{1/4} u||^2_{L^2(B_{2^{k-n-m}})}\nonumber \,.
\end{eqnarray}
We choose $\gamma ,\eta>0$ so that $C_1^{-1}C\gamma<1/4$ and $C_1^{-1}{\eta^2}<1/4\,.$ \par
We get for some constant $C>0$
\begin{eqnarray}\label{estrot7bis}
 &&|| M\Delta^{1/4} u||^2_{L^2(B_{2^{k-n-m}})}-\frac{1}{4}||\Delta^{1/4} u||^2_{L^2(B_{2^{k-n-m}})}\\[5mm]
 &&\le C\left[\sum_{h=k-n-m}^\infty (2^\frac{k-h}{2})||M\Delta^{1/4} u||^2_{L^2(A_h)}+  \sum_{h=k-n-m}^{+\infty}2^{\frac{k-h}{2}}|| \Delta^{1/4} u||_{L^2(A_h)}\right] \nonumber
\end{eqnarray}
Thus by replacing in \rec{estrot7bis} $k-n-m$ by $k$   we get \rec{estrot5} and we conclude the proof\,.
  \hfill$\Box$\par
 \medskip
 The {\bf proof of Theorem \ref{regularity2}} is analogous to that of  Theorem \ref{regularity1} and thus we omit it.


 \section{Morrey estimates and H\"older continuity of ${\mathbf {1/2}}$-Harmonic Maps into the Sphere}\label{application}
 
 We consider the $m-1$-dimensional sphere $S^{m-1}\subset\R^m$.   Let  $\Pi_{S^{m-1}}$ be the 
orthogonal projection on $S^{m-1}\,.$  
 We also consider the   {Dirichlet  energy}
 \begin{equation}\label{energy}
{\cal{L}}(u)= \int_{\R} |{\Delta}^{1/4} u(x)|^2 dx\,.
\end{equation} 
where $u\colon\R\to S^{m-1}$\,.\par 
The weak ${1/2}$-harmonic maps are defined as critical points
of the functional \rec{energy} with respect to perturbation of the form $\Pi_{S^{m-1}}(u+t\phi)$, where $\phi$ is
an arbitrary compacted supported smooth map    from $\R $ into $\R^m\,.$
\begin{Definition}\label{weakhalfharm}
We say that $u\in H^{1/2}(\R,S^{m-1})$ is a weak  ${1/2}$-harmonic map if and only if, for every
arbitrary compacted supported smooth maps $\phi$  from $\R $ into $\R^m$  we have
$$
\frac{d}{dt}{\cal{L}}(\Pi_{S^{m-1}}(u+t\phi))_{|t=0}\,.
$$
\end{Definition}
We introduce some notations.
 We denote by $\bigwedge(\R^m)$ the exterior algebra (or Grassmann Algebra) of $\R^m$ and by the symbol $\wedge$ the {\em exterior or wedge product}. 
For every $p=1,\ldots,m$,   $\bigwedge_p(\R^m)$ is the vector space of $p$-vectors  \par If $(e_i)_{i=1,\ldots,m}$ is the 
canonical orthonormal  basis of $\R^m$ then every element $v\in \bigwedge_p(\R^m)$ is written as
$v=\sum_{I}v_{I}e_{I}$ where $I=\{i_1,\ldots,i_p\}$ with $1\le i_1\le\ldots\le i_p\le m$ , $v_I:=v_{i_1,\ldots,i_p} $ and $ e_{I}=:=e_{i_1}\wedge\ldots\wedge e_{i_p}\,.$ \par
By the symbol $\res$ we denote the interior multiplication  $\res\colon \bigwedge_p(\R^m)\times \bigwedge_q(\R^m)\to\bigwedge_{q-p}(\R^m)$ defined as follows. \par
 Let   $e_I=e_{i_1}\wedge\ldots\wedge e_{i_p}$, $e_J=e_{j_1}\wedge\ldots\wedge e_{j_q},$ with $q\ge p\,.$ Then  $e_I\res e_J=0$ if $I\not\subset J$, otherwise
  $e_I\res e_J=(-1)^Me_{K}$ where $e_{K} $ is a $q-p$ vector and $M$ is the number of pairs $(i,j)\in I\times J$ with $j>i\,.$ \par 
By the symbol $\bullet$ we denote the first order contraction between multivectors. 
We recall that it satisfies $\alpha\bullet \beta =\alpha\res \beta$ if $\beta$ is a $1$-vector and 
 $\alpha\bullet( \beta\wedge\gamma) =(\alpha\bullet \beta)\wedge \gamma+(-1)^{pq}(\alpha\bullet\gamma)\wedge\beta $,  if $\beta$ and $\gamma$ are respectively a $p$-vector and a $q$-vector\,.
\par
Finally by the symbol $\ast$ we denote the Hodge-star operator, $\ast\colon\bigwedge_p(\R^m)\to\bigwedge_{m-p}(\R^m)$,
defined by $\ast\beta=(e_{1}\wedge\ldots\wedge e_n)\bullet \beta$.  
For an introduction of the Grassmann Algebra we refer the reader to the first Chapter of the book by Federer\cite{fed}\,.\par
 
Next we    write the Euler equation associated to the functional \rec{energy}\,.\par
\begin{Theorem}\label{eulerth}
All weak $1/2$-harmonic maps $u\in H^{1/2}(\R,{\cal{N}})$ satisfy in a weak sense
\par
i)  the equation
\begin{equation}\label{perp}
\int_{\R}{\Delta}^{1/2} u v \,dx=0,
\end{equation}
for every $v\in C^{\infty}_0(\R,\R^m)$ and $v\in T_{u(x)}{\cal{N}}$ almost everywhere, or in a equivalent way\par
ii) the equation 
\begin{equation}\label{wedge}
 {\Delta}^{1/2}u\wedge u=0,~~\mbox{in ${\cal{D}}^\prime\,,$}
 \end{equation}
 or 
\par
iii) the equation
\begin{equation}\label{euler1}
\Delta^{1/4}(u\wedge \Delta^{1/4} u)=T(Q ,u)\,;
\end{equation}
with $Q=u\wedge\,.$  
 \end{Theorem}
 {\bf\noindent Proof of Theorem \ref{eulerth}}\par
 The proof of \rec{perp} is the same of the one of Lemma 1.4.10 in \cite{Hel} , so we omit it.\par
  We prove \rec{wedge}. 
 We take $\varphi\in C_0^\infty(\R,\bigwedge_{m-2}(\R^m)).$  The following holds
 \begin{equation}\label{wedge2}
 \int_{\R} \varphi\wedge u\wedge {\Delta}^{1/2}u \,dx=\left(\int_{\R}  \ast(\varphi\wedge u)\cdot{\Delta}^{1/2}u\,dx\right)e_1\wedge\ldots\wedge e_m \end{equation}
 {\bf Claim~:}  $v= \ast(\varphi\wedge u)\in \dot {H}^{1/2}(\R,\R^m)$ and $v(x)\in T_{u(x)}S^{m-1}$ a.e.\par
 {\bf Proof of the claim.}\par
 The fact that $v\in \dot{H}^{1/2}(\R,\R^m)$ follows form the fact that its components are the product of two functions which are in
 $\dot{H}^{1/2}(\R,\R^m)\cap L^\infty(\R,\R^m)$\,.\par      
 We have
 \begin{equation}
   v\cdot u =  \ast (u\wedge\varphi)\cdot u=\ast(u\wedge\varphi\wedge u)=0\,.
 \end{equation}
 
 It follows from \rec{perp} and \rec{wedge2} that
 $$ \int_{\R} \varphi\wedge u\wedge {\Delta}^{1/2}u dx=0\,.$$
 This shows that  $ ({\Delta})^{1/2}u\wedge u=0,$ in ${\cal{D}}^\prime\,,$ and we can conclude\,.\par
  As far as equation \rec{euler1}  is concerned it is enough to observe that $\Delta^{1/2} u\wedge   u=0$ and  
  $\Delta^{1/4}u\wedge\Delta^{1/4} u=0\,.$ ~~\hfill$\Box$
  \par
  \medskip
   By combing Theorem \ref{eulerth},   Proposition 
  \ref{pr-I.2} and the results of the previous Section we get the H\"older regularity
   of weak $1/2$-harmonic maps.\par
   \begin{Theorem}\label{reghm}
Let  $u\in \dot{H}^{1/2}(\R,S^{m-1})$ be a harmonic map. Then  $u\in C^{0,\alpha}(\R,{S^{m-1}})\,.$
\end{Theorem}
  {\bf Proof of \ref{reghm}.} From Theorem \ref{eulerth} it follows that $u$ satisfies equation \rec{euler1}. Moreover, since $|u|=1$, Proposition 
  \ref{pr-I.2} implies  that $u$ satisfies  \rec{eq2intr} as well\,.
  Theorem \ref{regularity1} and Theorem \ref{regularity2} yield  respectively that for $k<0$, with $|k|$  large enough
   
   \begin{eqnarray}\label{estrot8}
&&
||u\wedge\Delta^{1/4} u||^2_{L^2 (B_{2^k})}\le C\sum_{h=k}^\infty (2^\frac{k-h}{2})||\Delta^{1/4} u||^2_{L^2(A_h)}+\frac{1}{4}||\Delta^{1/4} u||^2_{L^2 (B_{2^k})}\,.\ \end{eqnarray}
 and
 \begin{eqnarray}\label{estrot9}
&&
||u\cdot\Delta^{1/4} u||^2_{L^2 (B_{2^k})}\le C\sum_{h=k}^\infty (2^\frac{k-h}{2})||\Delta^{1/4} u||^2_{L^2(A_h)} +\frac{1}{4}||\Delta^{1/4} u||^2_{L^2 (B_{2^k})}\,.\end{eqnarray}
  Since
  $$|| \Delta^{1/4} u||^2_{L^2 (B_{2^k})}=||u\cdot\Delta^{1/4} u||^2_{L^2 (B_{2^k})}+||u\wedge\Delta^{1/4} u||^2_{L^2 (B_{2^k})} \,,$$
  we get
  \begin{eqnarray}\label{estrot10}|| \Delta^{1/4} u||^2_{L^2 (B_{2^k})}\le C\sum_{h=k}^\infty (2^\frac{k-h}{2})||\Delta^{1/4} u||^2_{L^2(A_h)}\,.\end{eqnarray}
  Now observe that for some  $C>0$ (independent on $k$) we have
  \begin{eqnarray}\label{locl2}
  &&
  C^{-1}\sum_{h=-\infty}^{k-1} ||\Delta^{1/4} u||^2_{L^2(A_h)}\le || \Delta^{1/4} u||^2_{L^2 (B_{2^k})}\le C \sum_{h=-\infty}^{k} ||\Delta^{1/4} u||^2_{L^2(A_h)}\,.\nonumber
  \end{eqnarray}
  Thus from \rec{locl2} and \rec{estrot9} it follows
  $$
  \sum_{h=-\infty}^{k-1} ||\Delta^{1/4} u||^2_{L^2(A_h)}\le C \sum_{h=k}^\infty (2^\frac{k-h}{2})||\Delta^{1/4} u||^2_{L^2(A_h)} \,.$$
  By applying Proposition \ref{seq} and using again \rec{estrot10} we get for $r>0$ small enough and some $\beta\in (0,1)$   
 
\begin{equation}\label{holdcond2}
\int_{B_r}|\Delta^{1/4}u|^2 dx\le Cr^\beta\,.\end{equation}
Condition \rec{holdcond2} yields that $u$ belongs to the Morrey-Campanato Space ${\cal{L}}^{2,-\beta}$ (see \cite{AD}), 
and thus $u\in C^{0,\beta/2}(\R)\,,$  (see for instance  \cite{AD, Gia})\,.~~\hfill$\Box$

   

\end{document}